\documentclass[10pt]{amsart}
\usepackage{latexsym}
\usepackage{amsfonts}
\usepackage{graphicx}
\usepackage{color}
\usepackage[active]{srcltx} 
\usepackage{booktabs} 
\usepackage{tikz}
\usepackage{amssymb}
\usepackage[font=small,labelfont=bf]{caption}
\usepackage{soul}\usepackage{geometry} 

\allowdisplaybreaks

\newcommand{\ext}{\raise1pt\hbox{$\textstyle\bigwedge$}}

\title{A geometric description of Discrete Exterior Calculus for general triangulations}

\author{Humberto Esqueda\and Rafael Herrera \and Salvador Botello\and Miguel Angel Moreles
}
\address{\newline Centro de Investigaci\'{o}n en Matem\'{a}ticas\newline
Jalisco s/n, Valenciana\newline
Guanajuato, GTO 36240,Mexico\newline
\textit{email: }esqueda,rherrera,botello,moreles@cimat.mx
}

\date{}

\vspace{-20pt}

\begin{document}

\begin{abstract}
We revisit the theory of Discrete Exterior Calculus (DEC) in 2D for general triangulations, relying only on Vector Calculus and Matrix Algebra. 
We present DEC numerical solutions of the Poisson equation and compare them against those 
found using the Finite Element Method with linear elements (FEML).
\end{abstract}

\maketitle

\tableofcontents

\section{Introduction}

The purpose of this paper is to introduce the theory of Discrete Exterior Calculus (DEC) 
to the widest possible audience and, therefore, we will 
rely mainly on Vector Calculus and Matrix Algebra.
Discrete Exterior Calculus is a relatively new method for solving partial differential equations \cite{HiraniThesis} based on the idea
of discretizing the mathematical theory of Exterior Differential Calculus, a theory that goes back to E. Cartan \cite{Cartan} and is fundamental 
in the areas of Differential Geometry and Differential Topology.
Although Exterior Differential Calculus is an abstract mathematical theory, it has been introduced in various fields such as in digital geometry processing \cite{Craneetal}, 
numerical schemes for partial differential equations \cite{HiraniThesis, Arnold}, etc.

In his PhD thesis \cite{HiraniThesis}, Hirani laid down the fundamental concepts 
of Discrete Exterior Calculus (DEC), using discrete combinatorial and geometric 
operations on simplicial complexes (in any dimension), proposing
discrete equivalents for differential forms, vector fields, differential and geometric operators, etc.
Perhaps the first numerical application of DEC to PDE was given in \cite{Hirani_K_N} in order to solve Darcy flow and Poisson's equation. 
In \cite{Griebel_R_S}, the authors develop a modification of DEC and show that in simple cases (e.g. flat geometry and regular meshes), the equations
resulting from DEC are equivalent to classical numerical schemes such as finite difference or finite volume
discretizations.
In \cite{Mohamedetal}, the authors used DEC to solve the Navier-Stokes equations and, 
in \cite{Dassiosetal} DEC was used with a discrete lattice model to simulate  elasticity, plasticity and failure of isotropic materials.

In this expository paper, we review the various
operators of Exterior Differential Calculus in 2D in terms of ordinary vector calculus, and introduce only the geometrical ideas that are essential 
to the formulation. 
Among those ideas is that of duality between the differentiation operator (on vector fields) and the boundary operator (on the domain) contained in Green's theorem. 
This duality is one of the key ideas of the method, which justifies taking the discretized derivative matrix as the transpose of the boundary operator matrix on the given mesh.  
Another important ingredient is the Hodge star operator, which is hidden in the notation of Vector Calculus. In order to show the necessity of the Hodge star operator, 
we carry out some simple calculations. In particular, we will introduce the notion of wedge product of vectors which, roughly speaking, helps us assign algebraic objects to parallelograms and
carry out algebraic manipulations with them. We present  DEC in the simplest terms possible using easy examples.
We also review the formulation of DEC for arbitrary meshes, which was first considered in \cite{Hirani2}.
Performance of the method is tested on the Poisson equation and compared with the Finite Element Method with linear elements (FEML).

The paper is organized as follows. 
In Section \ref{sec: preliminaries}, we introduce the wedge product of vectors and the {\em geometric} Hodge star operator, and
rewrite Green's theorem appropriately in order to display the duality between the differentiation and
the boundary operators.
In Section \ref{sec: discrete exterior calculus}, we present the operators of DEC (mesh, dual mesh, discrete derivation, discrete Hodge star operator), showing simple examples throughout.
In Section \ref{sec: general triangle}, we present  the formulation of DEC on arbitrary triangulations.
In Section \ref{sec: example}, we present the numerical solution of a Poisson equation with DEC and FEML, in order to compare their performance.
In Section \ref{sec: conclusions}, we present our conclusions.

\section{2D Exterior Differential Calculus as Vector Calculus}
\label{sec: preliminaries}

In this section we introduce two geometric operators (the wedge product and the Hodge star) and explain how to use them together with the gradient operator in order to obtain the Laplacian.

\subsection{Wedge product for vectors in $\mathbb{R}^2$}

Let $a,b$ be vectors in $\mathbb{R}^2$. We can assign to these vector the parallelogram
they span, and to such a parallelogram its area. The latter is equal to the determinant of the transformation matrix  sending $e_1,e_2$ to $a,b$ respectively. In Exterior Calculus, such a parallelogram is regarded as an algebraic object $a\wedge b$, a bivector, and the set of bivectors is equipped with a vector space structure to form the one-dimensional vector space 
\[\ext^2\mathbb{R}^2.\]
Thus, we have the following spaces
\[\mathbb{R}, \quad\quad\mathbb{R}^2, \quad\quad\ext^2\mathbb{R}^2,\]
of scalars, vectors and bivectors, respectively.
For instance, the square formed by $e_1$ and $e_2$ is represented by the symbol
\[e_1\wedge e_2\]
which is read as "$e_1$ wedge $e_2$" (see Figure \ref{fig: square})
    \begin{center}
    \includegraphics[width=0.2\textwidth]{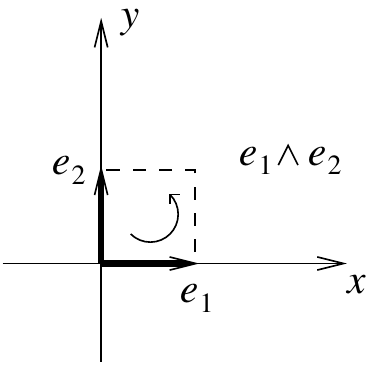}%
    \captionof{figure}{Wedge product of two vectors.}
    \label{fig: square}
    \end{center}
In $\mathbb{R}^2$, this represents an ``element'' of unit area.

Note that if we list the vectors in the opposite order, we have a different orientation and, therefore, the algebraic objects must satisfy $e_2\wedge e_1=-e_1\wedge e_2$ (see Figure \ref{fig: -square})
    \begin{center}
    \includegraphics[width=0.26\textwidth]{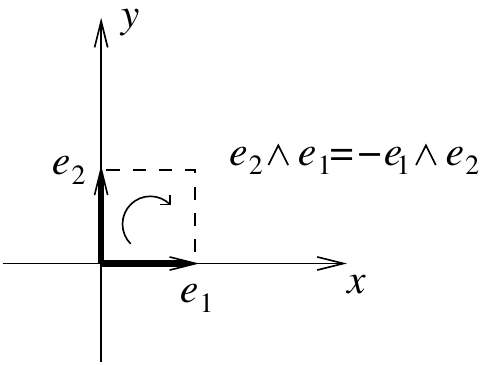}%
    \captionof{figure}{Change of orientation of the parallelogram implies anticommutativity in the wedge product of two vectors.}\label{fig: -square}
    \end{center}
More generally, given two vectors $a,b\in\mathbb{R}^2$, their wedge product $a\wedge b$ looks as follows (see Figure \ref{fig: a wedge b})
\begin{center}
 \includegraphics[width=0.25\textwidth]{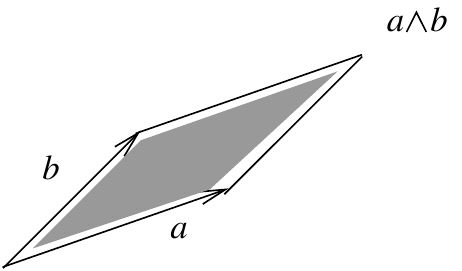}
 \captionof{figure}{The wedge product of the vectors $a$ and $b$}\label{fig: a wedge b}
\end{center}
The properties of the wedge product are
\begin{itemize} 
\item it is anticommutative: $a\wedge b = - b\wedge a$ (see Figure \ref{fig: anticommutativity}
\begin{center}
 \includegraphics[width=0.25\textwidth]{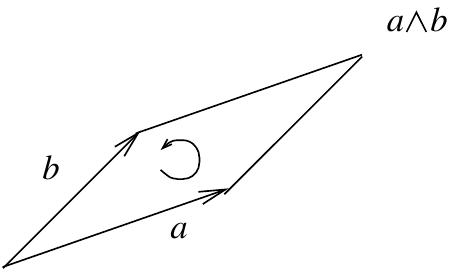}
 \hspace{.5in}
 \includegraphics[width=0.25\textwidth]{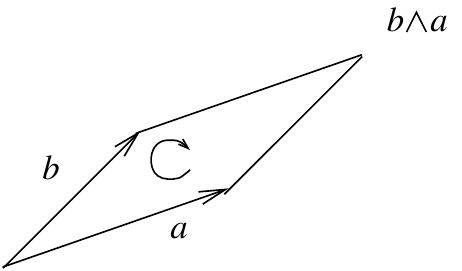}
 \captionof{figure}{Anticommutativity of the wedge product.}\label{fig: anticommutativity}
 \end{center}
\item $a\wedge a = 0$ since it is a parallelogram with area zero (see Figure \ref{fig: zero divisor})
\begin{center}
 \includegraphics[width=0.25\textwidth]{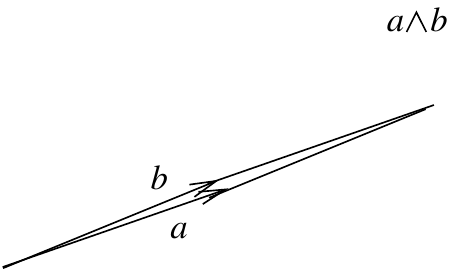}
 \captionof{figure}{Parallelogram with very small area, depicting what happens when $b$ tends to $a$.}\label{fig: zero divisor}
\end{center}
 \item it is distributive
\[(a + b )\wedge c = a\wedge c + b\wedge c;\]
 \item it is associative
\[(a \wedge b )\wedge c = a\wedge (b \wedge  c).\]
\end{itemize}

For example, let
\begin{eqnarray*}
 a&=& (a_1,a_2)\,\, = a_1 e_1 + a_2 e_2,\\
 b&=& (b_1,b_2)\,\, \,= b_1 e_1 + b_2 e_2.
\end{eqnarray*}
Then
\begin{eqnarray*}
 a\wedge b 
    &=&
   (a_1 e_1 + a_2 e_2)\wedge (b_1 e_1 + b_2 e_2)\\
    &=&
   a_1b_1\, e_1\wedge e_1 + a_1b_2\, e_1\wedge e_2 + a_2b_1\, e_2\wedge e_1 + a_2b_2\, e_2\wedge e_2 \\
    &=&
    a_1b_2\, e_1\wedge e_2 - a_2b_1\, e_1\wedge e_2  \\
    &=&
    (a_1b_2 - a_2b_1)\, e_1\wedge e_2 ,
\end{eqnarray*}
i.e. the determinant 
\[\det\left(\begin{array}{cc}
a_1 & a_2 \\ 
b_1 & b_2
\end{array} \right)\]
times the canonical bivector/square $e_1\wedge e_2$.

\subsection{Hodge star operator for vectors in $\mathbb{R}^2$}\label{subsec: Hodge star}
Now consider the following situation: given the vector $e_1$, find another vector $v$ such that the parallelogram that they form has area $1$. It is readily seen that, for instance, $v=e_2, -e_2,e_1+e_2$ are all solutions. Requiring orthogonality and standard orientation, we see that 
$e_2$ is the unique solution.
This process is summarized in the Hodge star operator, which basically says that the complementary vector for $e_1$ is $e_2$, and the one for $e_2$ is $-e_1$,
\begin{eqnarray*}
\star e_1 &=& e_2\\
{}\star e_2 &=& -e_1.
\end{eqnarray*}

In general, the equation that defines the Hodge star operator for any given vector $v\in \mathbb{R}^2$ is the following
\[w\wedge (\star v) = (w\cdot v) \, e_1\wedge e_2.\]
for every $w\in \mathbb{R}^2$.
In particular, if we take $v=w$,
\[v\wedge (\star v) = |v|^2 \, e_1\wedge e_2,\]
which means that $v$ and $\star v$ form a square of area $|v|^2$.
Thus, the Hodge star operator on a vector $v\in \mathbb{R}^2$ can be thought of as finding the vector that makes with $v$ a positively oriented square of area $\vert v\vert^2$.

It is somewhat less intuitive to work out the Hodge star of a bivector.
First of all, we have to treat bivectors as vectors of a different space, namely the space of bivectors $\ext^2\mathbb{R}^2$.
Secondly, the length of a bivector $v\wedge w$ is its area
\[{\rm length}(v\wedge w):= {\rm Area}(v\wedge w).\]
Thus, the defining equation of the Hodge star applied to $(v\wedge w)$ and  $\star(v\wedge w)$ reads as follows
\begin{eqnarray*}
(v\wedge w)\wedge \star(v\wedge w)&=& \left<v\wedge w, v\wedge w\right>e_1\wedge e_2
\end{eqnarray*}
which means
\begin{eqnarray*}
(v\wedge w)\wedge \star(v\wedge w)
&=& {\rm Area}(v\wedge w)^2 e_1\wedge e_2.
\end{eqnarray*}
Since $v\wedge w={\rm Area}(v\wedge w) e_1\wedge e_2$ is already a bivector,
$\star(v\wedge w)$ must be a scalar, i.e.
\[\star(v\wedge w)={\rm Area}(v\wedge w).\]
When $v=e_1$ and $w=e_2$ we have
\[\star(e_1\wedge e_2) = 1.\]

Finally, the Hodge star of a number $\lambda$ is a bivector, i.e.
\[\star\lambda=\lambda \, e_1\wedge e_2.\]

\subsection{The Laplacian}
Let $f:\mathbb{R}^2\rightarrow \mathbb{R}$ and consider the gradient
\[\nabla f = {\partial f\over \partial x}e_1+{\partial f\over \partial y}e_2.\] 
Apply the Hodge star operator to it
\begin{eqnarray*}
\star\nabla f 
   &=& {\partial f\over \partial x}*e_1+{\partial f\over \partial y}*e_2\\
   &=& {\partial f\over \partial x}e_2-{\partial f\over \partial y}e_1.
\end{eqnarray*} 
Now take the gradient of each coefficient function together with wedge product
\begin{eqnarray*}
\nabla\wedge\star\nabla f 
   &:=& \nabla\left({\partial f\over \partial x}\right)\wedge e_2-\nabla\left({\partial f\over \partial y}\right)\wedge e_1\\
   &=& \left({\partial^2 f\over \partial x^2}e_1+{\partial^2 f\over \partial y\partial x}e_2\right)\wedge e_2
   -\left({\partial^2 f\over \partial x\partial y}e_1+{\partial^2  f\over \partial^2 y}e_2\right)\wedge e_1\\
   &=& {\partial^2 f\over \partial x^2}e_1\wedge e_2+{\partial^2 f\over \partial y\partial x}e_2\wedge e_2
   -{\partial^2 f\over \partial x\partial y}e_1\wedge e_1-{\partial^2  f\over \partial^2 y}e_2\wedge e_1\\
   &=& \left({\partial^2 f\over \partial x^2}+{\partial^2  f\over \partial^2 y}\right)e_1\wedge e_2.
\end{eqnarray*} 
By taking the Hodge star of this last expression we get the ordinary Laplacian of $f$
\begin{eqnarray*}
\star\nabla\wedge\star\nabla(f) 
   &=& {\partial^2 f\over \partial x^2}+{\partial^2  f\over \partial^2 y}.
\end{eqnarray*}

\noindent{\bf Remarks. } \newline
(i) This rather convoluted looking way of computing the Laplacian of a function is based on Exterior Differential Calculus, a theory that generalizes 
the operators of vector calculus (gradient, curl and divergence) to arbitrary dimensions, and is the basis for the differential topological theory of deRham cohomology.\newline
(ii) We would like to emphasize the  necessity of using the Hodge star operator  $*$ in order to make the combination of differentiation 
and wedge product produce the correct answer.

\subsection{Duality in Green's theorem}

 Green's theorem states that for a vector field $(L,M)$ defined on a region $D\subset  \mathbb{R}^2$,
 \[\int_D\left({\partial M\over \partial x}-{\partial L\over \partial y}\right)dxdy=\int_{C=\partial D} (Ldx+Mdy) .\]
In this section, we will explain how this identity encodes a duality between the operator of differentiation
and that of taking the boundary of the domain of integration.

\subsubsection{Rewriting Green's theorem}

Note that if $F=(L,M)$ is a vector field, we can write it as
\[F= Le_1+Me_2.\]
Then, we can apply the gradient operator together with wedge product in the following fashion
\begin{eqnarray*}
\nabla^\wedge (L e_1 +  M  e_2)
   &:=&
\nabla L\wedge e_1 + \nabla M \wedge e_2\\
   &=&
  \left({\partial L\over \partial x}e_1+{\partial L\over \partial y}e_2\right)\wedge e_1
  +\left({\partial M\over \partial x}e_1+{\partial M\over \partial y}e_2\right)\wedge e_2\\
   &=&
  {\partial L\over \partial x}e_1\wedge e_1+{\partial L\over \partial y}e_2\wedge e_1
  +{\partial M\over \partial x}e_1\wedge e_2+{\partial M\over \partial y}e_2\wedge e_2\\   
   &=&
  {\partial L\over \partial y}e_2\wedge e_1
  +{\partial M\over \partial x}e_1\wedge e_2\\   
   &=&
  \left({\partial M\over \partial x}
  -{\partial L\over \partial y}\right)e_1\wedge e_2.
\end{eqnarray*}
Note that we have defined a new operator $\nabla^\wedge$ which combines differentiation and wedge product.
Applying the Hodge star operator we obtain
\begin{eqnarray*}
\star\nabla^\wedge (L e_1 +  M e_2)
   &=&
  \left({\partial M\over \partial x}
  -{\partial L\over \partial y}\right),
\end{eqnarray*}
i.e. the integrand of the left hand side of the identity of integrals in Green's theorem.
Thus, as a first step, Green's theorem can be rewritten as follows:
\[\int_D\nabla^\wedge (L, M)dxdy=\int_{\partial D} (L ,  M )\cdot (dx,dy)\]

\subsubsection{The duality in Green's theorem}

 Let us recall the following fact from Linear Algebra. Given a linear transformation $A$ of Euclidean space $\mathbb{R}^n$,
the transpose $A^T$ satisfies 
\[\left<A(v), w \right>= \left<v, A^T(w)\right>\]
for any two vectors $v,w\in\mathbb{R}^n$, where $\left<\cdot, \cdot\right>$ denotes the standard inner/dot product. 
In fact, such an identity characterizes the transpose $A^T$.

Now, let us do the following notational trick: substitute the integration symbols in Green's theorem by $\ll\cdot,\cdot\gg$ as follows:
\begin{eqnarray*}
\int_C (L,M)\cdot(dx,dy) &=& \ll (L,M), C\gg,\\
\int_D\nabla^\wedge (L, M)dxdy &=& \ll\nabla^\wedge(L,M), D\gg,
\end{eqnarray*}
where $C=\partial D$ is the boundary of the region.
Using this notational change, Green's theorem reads as follows
\[\ll\nabla^\wedge(L,M), D\gg\,\,=\,\,\ll(L,M), \partial D\gg.\]
Roughly speaking, this means that the differential operator $\nabla^\wedge$ is the transpose of the boundary operator $\partial$ by means of the product $\ll\cdot,\cdot\gg$.

{\bf Remark}. The previous observation is fundamental in the development of DEC, since the boundary operator is well understood and easy to calculate on meshes.

\section{Discrete Exterior Calculus}\label{sec: discrete exterior calculus}

Now we will discretize the differentiation operator $\nabla^\wedge$ presented above. We will start by describing the discrete version of the boundary operator
on simplices/triangles. Afterwards, we will treat the differentiation operator as the transpose of the boundary operator. 

We are interested in using certain geometric subsets of a given triangular mesh of a 2D region.
Such subsets include vertices/nodes, edges/sides and faces/triangles.
We will describe each one of them by means of the ordered list of vertices whose convex closure constitutes the subset of interest.
For instance, consider the triangular mesh of the planar hexagonal region in Figure \ref{fig: hexagon01},
    \begin{center}
    \includegraphics[width=0.3\textwidth]{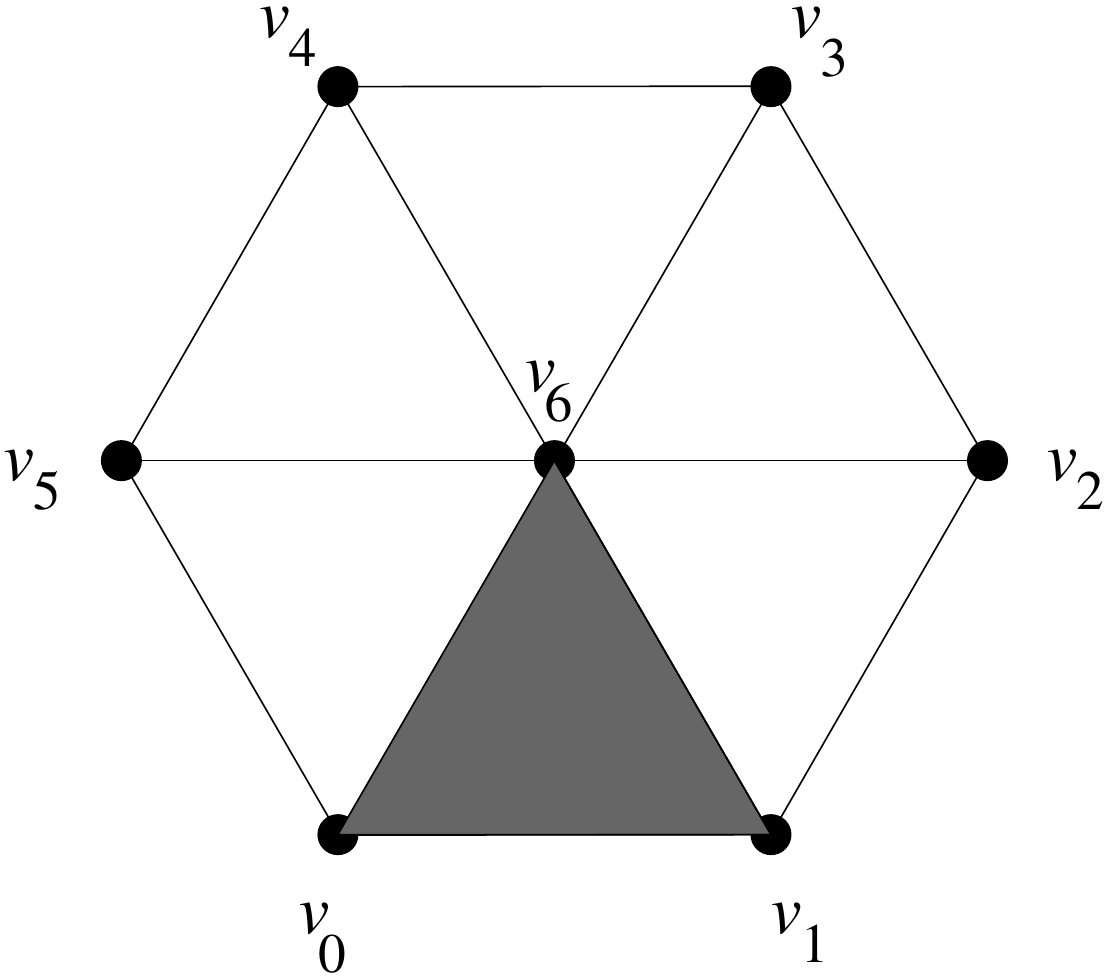}%
    \captionof{figure}{Triangular mesh of a planar hexagonal region.}\label{fig: hexagon01}
    \end{center}
where the shaded triangle will be denoted by $[v_0,v_1,v_6]$, and its edge joining the vertices $v_0$ and $v_1$ will be denoted by $[v_0,v_1]$.
For the sake of notational consistency, we will denote the vertices also enclosed in brackets, e.g. $[v_0]$.

\subsection{Boundary operator}

There is a well known boundary operator $\partial$ for oriented triangles, edges and points:
\begin{itemize} 
\item For points/vertices:
    \begin{center}
    \includegraphics[width=.07\textwidth]{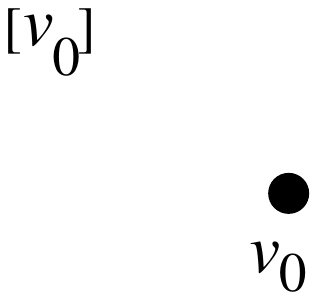}%
    \captionof{figure}{Boundary of a vertex: $\partial [v_0] = 0$.}
    \end{center}

\item For sides/edges:
    \begin{center}
    \includegraphics[width=0.2\textwidth]{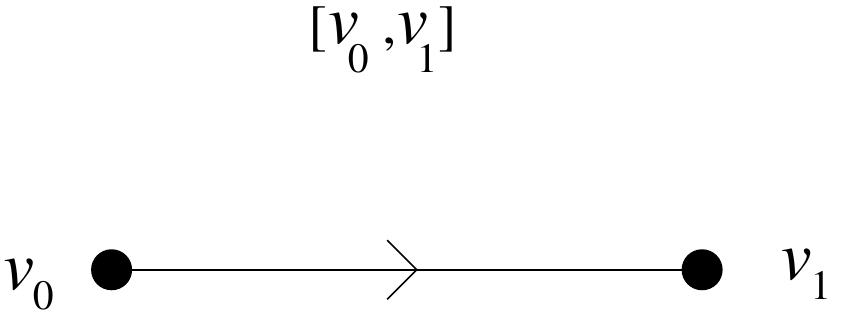}%
    \captionof{figure}{Boundary of an edge: $\partial [v_0,v_1] = [v_1] - [v_0]$.}
    \end{center}

\item For faces/triangles:
    \begin{center}
    \includegraphics[width=0.2\textwidth]{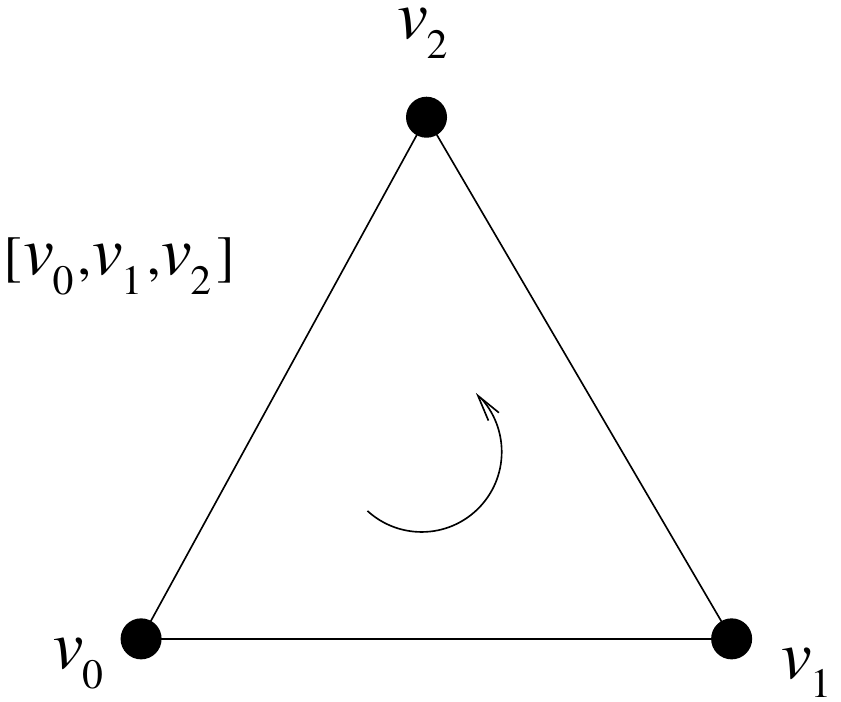}%
    \captionof{figure}{Boundary of a face: $\partial [v_0,v_1,v_2] = [v_1,v_2] - [v_0,v_2]+[v_0,v_1]$.}
    \end{center}

\end{itemize}

{\bf Example}.
Let us consider again the mesh of the planar hexagonal (with oriented triangles) in Figure \ref{fig: hexagon02}
    \begin{center}
    \includegraphics[width=0.3\textwidth]{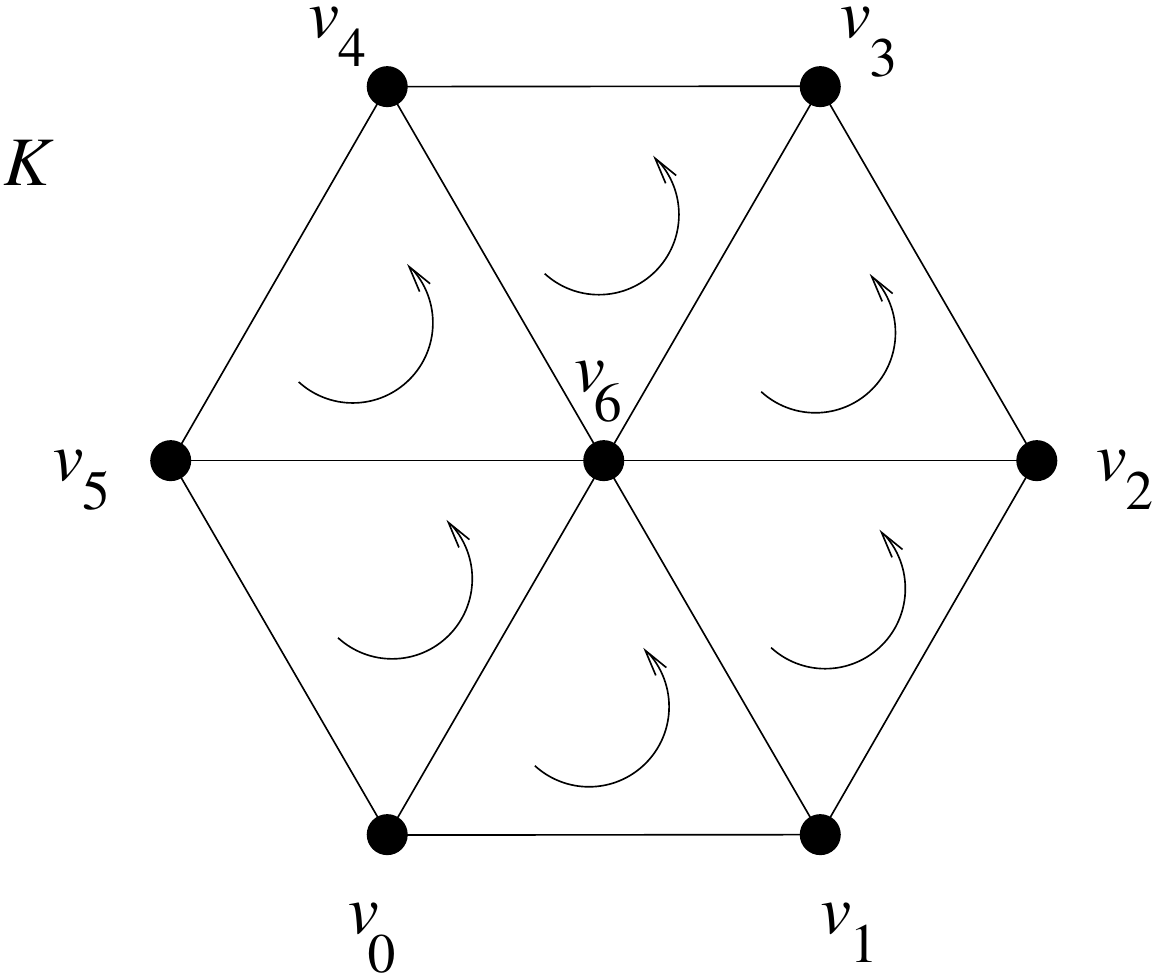}%
    \captionof{figure}{Oriented triangular mesh of a planar hexagonal region.}\label{fig: hexagon02}
    \end{center}
We will denote a triangle by the list of its vertices listed in order according to the orientation of the tringle.
Thus, we have the following ordered lists: 
\begin{itemize}
 \item list of faces
\[
\{[v_0,v_1,v_6],[v_1,v_2,v_6],[v_2,v_3,v_6],
[v_3,v_4,v_6],[v_4,v_5,v_6],[v_5,v_0,v_6]\};
\]    
\item list of edges
\[
\{
[v_0,v_6],
[v_1,v_6],
[v_2,v_6],
[v_3,v_6],
[v_4,v_6],
[v_5,v_6],
[v_0,v_1],
[v_1,v_2],
[v_2,v_3],
[v_3,v_4],
[v_4,v_5],
[v_5,v_0]
\};
\]
\item list of vertices
\[
\{
[v_0],[v_1],[v_2],[v_3],[v_4],[v_5],[v_6]
\}.
\]

\end{itemize}

A key idea in DEC is to consider each face as an element of a basis of a vector space. Namely,  coordinate vectors are associated to faces as follows:

\begin{eqnarray*}
[v_0,v_1,v_6] &\longleftrightarrow& (1,0,0,0,0,0),\\
{}[v_1,v_2,v_6] &\longleftrightarrow& (0,1,0,0,0,0),\\
{}[v_2,v_3,v_6] &\longleftrightarrow& (0,0,1,0,0,0),\\
{}[v_3,v_4,v_6] &\longleftrightarrow& (0,0,0,1,0,0),\\
{}[v_4,v_5,v_6] &\longleftrightarrow& (0,0,0,0,1,0),\\
{}[v_5,v_0,v_6] &\longleftrightarrow& (0,0,0,0,0,1).
\end{eqnarray*}
Similarly,  coordinate vectors are associated to the edges
\begin{eqnarray*}
{} [v_0,v_6] &\longleftrightarrow& (1,0,0,0,0,0,0,0,0,0,0,0),\\
{} [v_1,v_6] &\longleftrightarrow& (0,1,0,0,0,0,0,0,0,0,0,0),\\
{} [v_2,v_6] &\longleftrightarrow& (0,0,1,0,0,0,0,0,0,0,0,0),\\
{} [v_3,v_6] &\longleftrightarrow& (0,0,0,1,0,0,0,0,0,0,0,0),\\
{} [v_4,v_6] &\longleftrightarrow& (0,0,0,0,1,0,0,0,0,0,0,0),\\
{} [v_5,v_6] &\longleftrightarrow& (0,0,0,0,0,1,0,0,0,0,0,0),\\
{} [v_0,v_1] &\longleftrightarrow& (0,0,0,0,0,0,1,0,0,0,0,0),\\
{} [v_1,v_2] &\longleftrightarrow& (0,0,0,0,0,0,0,1,0,0,0,0),\\
{} [v_2,v_3] &\longleftrightarrow& (0,0,0,0,0,0,0,0,1,0,0,0),\\
{} [v_3,v_4] &\longleftrightarrow& (0,0,0,0,0,0,0,0,0,1,0,0),\\
{} [v_4,v_5] &\longleftrightarrow& (0,0,0,0,0,0,0,0,0,0,1,0),\\
{} [v_5,v_0] &\longleftrightarrow& (0,0,0,0,0,0,0,0,0,0,0,1).
\end{eqnarray*}
Finally, we do the same with the vertices
\begin{eqnarray*}
{} [v_0] &\longleftrightarrow& (1,0,0,0,0,0,0),\\
{} [v_1] &\longleftrightarrow& (0,1,0,0,0,0,0),\\
{} [v_2] &\longleftrightarrow& (0,0,1,0,0,0,0),\\
{} [v_3] &\longleftrightarrow& (0,0,0,1,0,0,0),\\
{} [v_4] &\longleftrightarrow& (0,0,0,0,1,0,0),\\
{} [v_5] &\longleftrightarrow& (0,0,0,0,0,1,0),\\
{} [v_6] &\longleftrightarrow& (0,0,0,0,0,0,1).
\end{eqnarray*}

Now, if we take the boundary of each face, we have
\begin{eqnarray*}
\partial [v_0,v_1,v_6] &=& [v_1,v_6] - [v_0,v_6]+[v_0,v_1],\\
\partial [v_1,v_2,v_6] &=& [v_2,v_6] - [v_1,v_6]+[v_1,v_2],\\
\partial [v_2,v_3,v_6] &=& [v_3,v_6] - [v_2,v_6]+[v_2,v_3],\\
\partial [v_3,v_4,v_6] &=& [v_4,v_6] - [v_3,v_6]+[v_3,v_4],\\
\partial [v_4,v_5,v_6] &=& [v_5,v_6] - [v_4,v_6]+[v_4,v_5],\\
\partial [v_5,v_0,v_6] &=& [v_0,v_6] - [v_5,v_6]+[v_5,v_0],
\end{eqnarray*}
which, under the previous assignments of coordinate vectors, corresponds to the linear transformation given by the following matrix
\[\partial_{2,1}=\left(\begin{array}{rrrrrr}
-1 & 0 & 0 & 0 & 0 & 1\\
1  & -1 & 0 & 0 &0 & 0\\
 0  & 1 & -1 & 0 & 0 & 0\\
 0  & 0 & 1 & -1 & 0 & 0\\
 0  &0  &0  & 1 & -1 & 0\\
 0  & 0 & 0 & 0 & 1 & -1\\
1  & 0 & 0 & 0 & 0 & 0\\
0   & 1 & 0 & 0 & 0 & 0\\
 0  & 0 & 1 &0  &0  &0 \\
  0 & 0 & 0 & 1 & 0 & 0\\
 0  & 0 & 0 & 0 & 1 & 0\\
 0  & 0 & 0 & 0 & 0 & 1
        \end{array}
\right),
\]
where the subindices in $\partial_{2,1}$ indicate that we are taking the boundary of 2-dimensional elements and obtaining 1-dimensional ones.
Similarly, taking the boundaries of all the edges gives
\begin{eqnarray*}
\partial [v_0,v_6] &=& [v_6] - [v_0],\\
\partial [v_1,v_6] &=& [v_6] - [v_1],\\
\partial [v_2,v_6] &=& [v_6] - [v_2],\\
\partial [v_3,v_6] &=& [v_6] - [v_3],\\
\partial [v_4,v_6] &=& [v_6] - [v_4],\\
\partial [v_5,v_6] &=& [v_6] - [v_5],\\
\partial [v_0,v_1] &=& [v_1] - [v_0],\\
\partial [v_1,v_2] &=& [v_2] - [v_1],\\
\partial [v_2,v_3] &=& [v_3] - [v_2],\\
\partial [v_3,v_4] &=& [v_4] - [v_3],\\
\partial [v_4,v_5] &=& [v_5] - [v_4],\\
\partial [v_5,v_0] &=& [v_0] - [v_5],
\end{eqnarray*}
which, under the previous assignments of coordinate vectors, corresponds to the linear transformation given by the following matrix
\[\partial_{1,0}=\left(
\begin{array}{rrrrrr|rrrrrr}
-1 & 0 & 0 & 0 & 0 & 0 & -1 & 0 & 0 & 0 & 0 & 1\\
0 & -1 & 0 & 0 & 0 & 0 & 1 & -1 & 0 & 0 & 0 & 0\\
0 & 0 & -1 & 0 & 0 & 0 & 0 & 1 & -1 & 0 & 0 & 0\\
0 & 0 & 0 & -1 & 0 & 0 & 0 & 0 & 1 & -1 & 0 & 0\\
0 & 0 & 0 & 0 & -1 & 0 & 0 & 0 & 0 & 1 & -1 & 0\\
0 & 0 & 0 & 0 & 0 & -1 & 0 & 0 & 0 & 0 & 1 & -1\\
1 & 1 & 1 & 1 & 1 & 1 & 0 & 0 & 0 & 0 & 0 & 0
\end{array}
\right).
\]

{\bf Remark}. Note how these matrices encode different levels of connectivity with orientations, such as who are the edges of which oriented triangle, or which are the 
end points of a given oriented edge. 

Due to the duality between $\partial$ and $\nabla^\wedge$, we can define the discretiztion of $\nabla^\wedge$ 
by
    \[\nabla^\wedge := (\partial)^T\]

For instance, we see that
the operator $\nabla^\wedge_{0,1}$ reads as follows
\[\nabla^\wedge_{0,1}=\left(\begin{array}{rrrrrrr}
-1 & 0 & 0 & 0 & 0 & 0 & 1\\
0 & -1 & 0 & 0 & 0 & 0 & 1\\
0 & 0 & -1 & 0 & 0 & 0 & 1\\
0 & 0 & 0 & -1 & 0 & 0 & 1\\
0 & 0 & 0 & 0 & -1 & 0 & 1\\
0 & 0 & 0 & 0 & 0 & -1 & 1\\
-1 & 1 & 0 & 0 & 0 & 0 & 0\\
0 & -1 & 1 & 0 & 0 & 0 & 0 \\
0 & 0 & -1 & 1 & 0 & 0 & 0 \\
0 & 0 & 0 & -1 & 1 & 0 & 0 \\
0 & 0 & 0 & 0 & -1 & 1& 0 \\
1 & 0 & 0 & 0 & 0 & -1& 0 
          \end{array}
\right).\]

\subsection{Dual mesh}

In order to discretize the Hodge star operator, we must first introduce  the notion of the dual mesh of a triangular mesh.

Consider the triangular mesh $K$ in Figure \ref{fig: dual mesh construction}(a).
The construction of the dual mesh $K^*$ is carried out as follows:
\begin{itemize}
 \item The vertices of the dual mesh $K^*$ are the circumcenters of the faces/triangles of the original mesh  (the blue dots in Figures 
 \ref{fig: dual mesh construction}(b) and \ref{fig: dual mesh construction}(c)). 
For instance, the dual of the face $[v_0,v_1,v_2]$ will be denoted by $[v_0,v_1,v_2]^*$.
 \item The edges of $K^*$ are the straight line segments joining the circumcenters of two adjacent triangles (those which share an edge).
Note that the resulting line segments 
 are orthogonal to one of the original edges (the blue straight line segments in 
 in Figures 
 \ref{fig: dual mesh construction}(b) and \ref{fig: dual mesh construction}(c)).
For instance, the dual of the edge $[v_0,v_6]$ will be denoted by $[v_0,v_6]^*$.
 \item The faces or cells of $K^*$ are the areas enclosed by the new polygons determined by the new edges.
For instance, the dual of the edge $[v_6]$ is the inner blue hexagon in 
in Figures 
 \ref{fig: dual mesh construction}(b) and \ref{fig: dual mesh construction}(c) 
and will be denoted by $[v_6]^*$.
\end{itemize} 
    \begin{center}
    {\small (a)}\includegraphics[width=0.28\textwidth]{boundary04-eps-converted-to.pdf}%
    \hspace{.25in}
    {\small (b)}\includegraphics[width=0.28\textwidth]{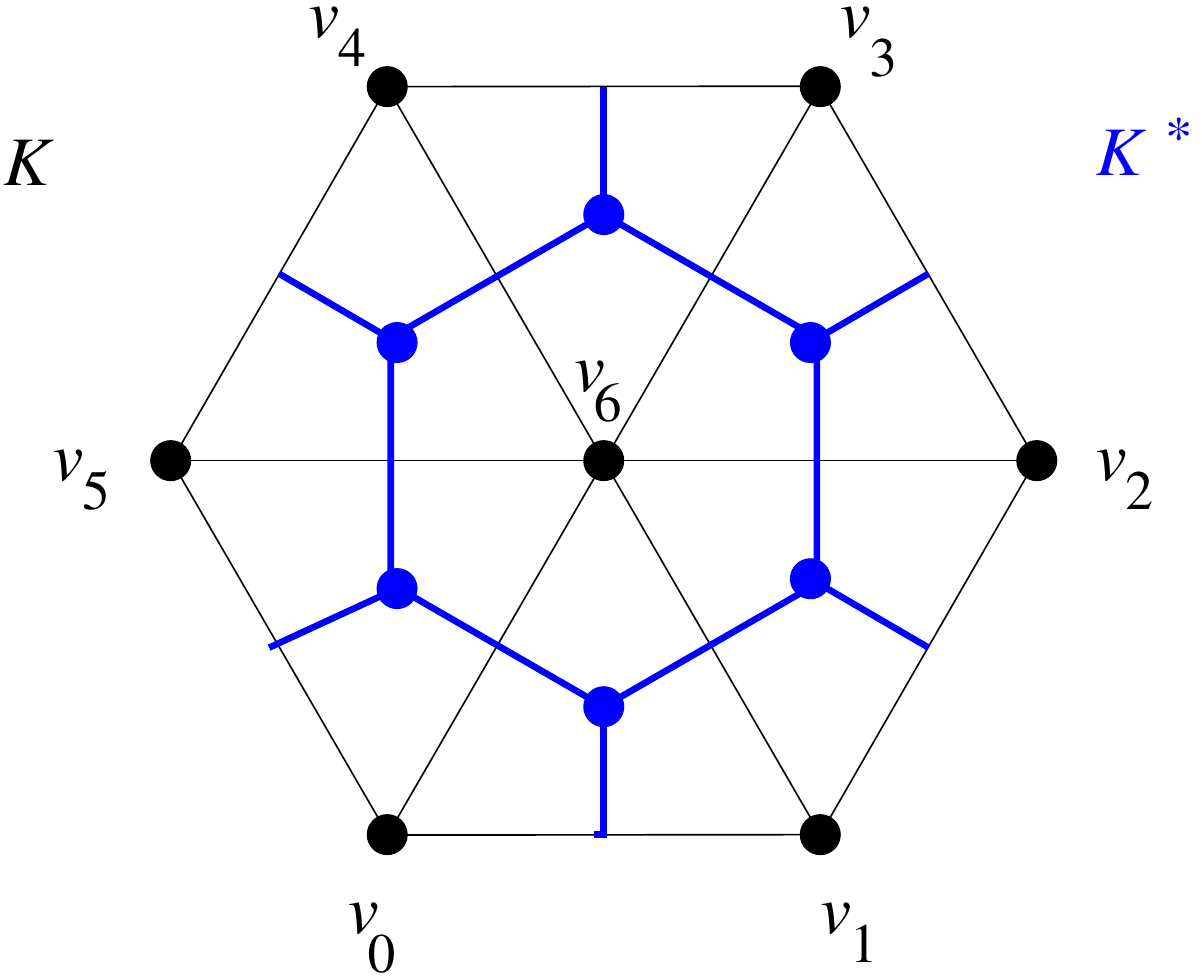}%
    \hspace{.25in}
    {\small (c)}\raise14pt\hbox{\includegraphics[width=0.22\textwidth]{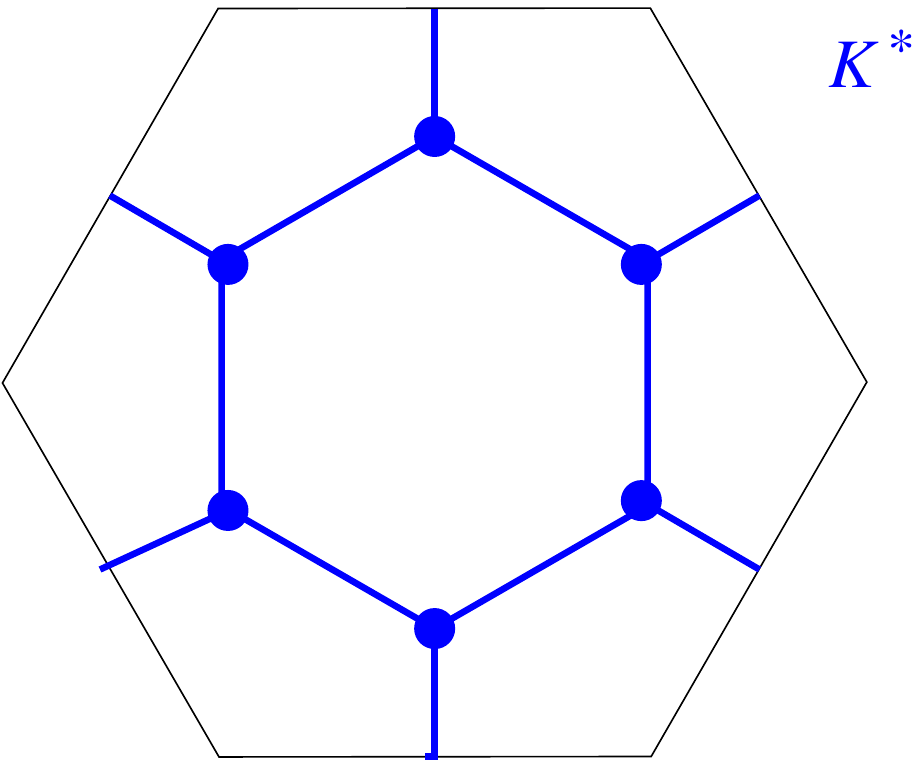}}%
    \captionof{figure}{Dual mesh construction: (a) triangular mesh $K$; (b) dual mesh $K^*$ superimposed on the mesh $K$; 
    (c) dual mesh $K^*$.}
    \label{fig: dual mesh construction}
    \end{center}

The orientation of the dual edges is given by the following recipe. If we have two adjacent triangles oriented as in 
Figure \ref{fig: two triangles 1}(a),
the dual edge crossing the edge of adjacency is oriented as in Figure \ref{fig: two triangles 1}(b)
    \begin{center}
    (a)\includegraphics[width=0.1\textwidth]{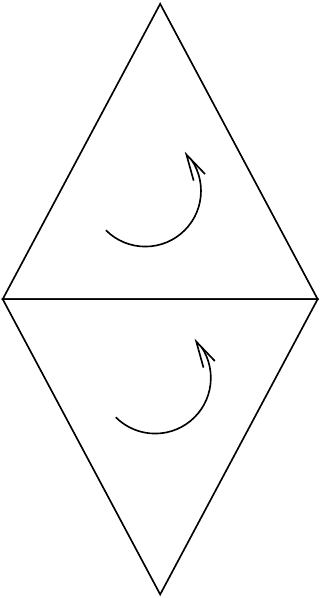}%
    \hspace{.8in}
    (b)\includegraphics[width=0.1\textwidth]{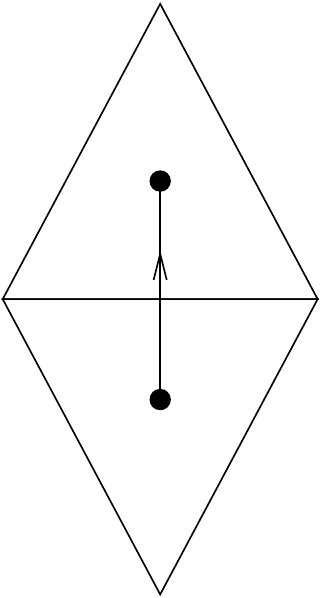}%
    \captionof{figure}{(a) Two adjacent oriented triangles. (b) Compatibly oriented dual edge.}\label{fig: two triangles 1}
    \end{center}

\subsection{Boundary operator on the dual mesh}

Consider the dual mesh in Figure \ref{fig: hexagon03} with the given labels and orientations
    \begin{center}
    \includegraphics[width=0.22\textwidth]{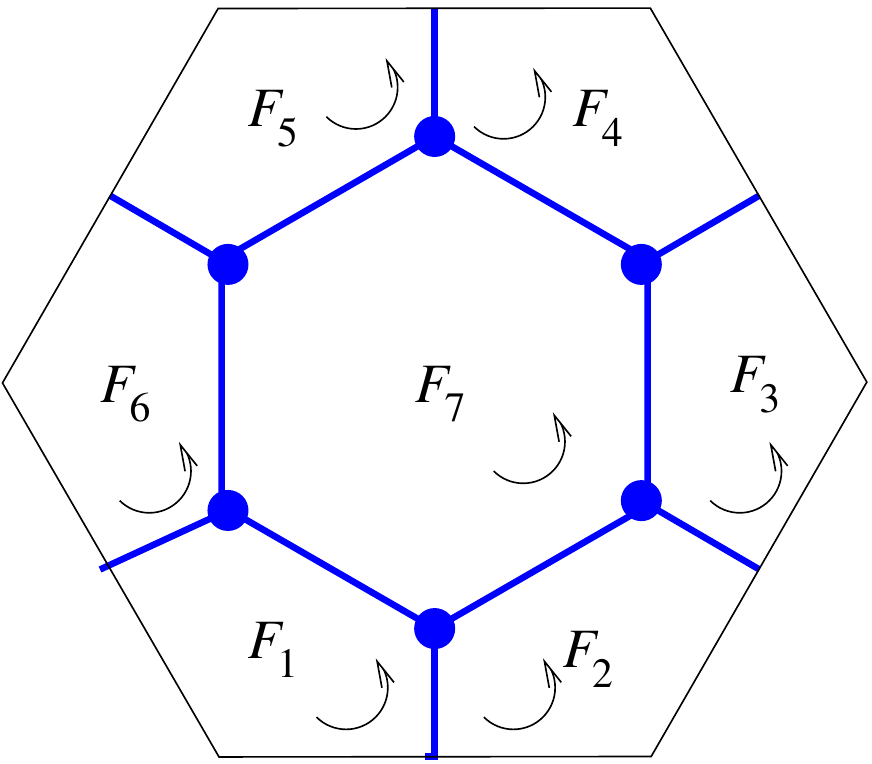}%
    \hspace{.5in}
    \includegraphics[width=0.22\textwidth]{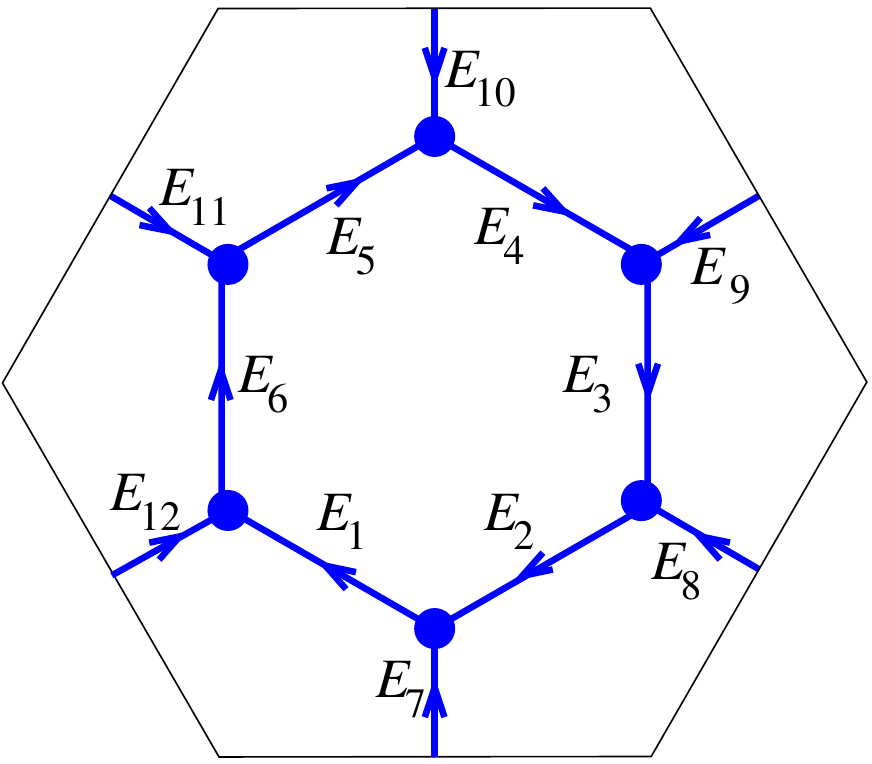}%
    \captionof{figure}{Oriented dual mesh.}\label{fig: hexagon03}
    \end{center}
The boundary operator is applied in a similar fashion as it was applied to triangles. In this case we have 
\begin{eqnarray*}
 \partial^{dual}_{2,1} F_1 &=&  E_{1}  + E_{7}  -  E_{12},\\
 \partial^{dual}_{2,1} F_2 &=&  E_{2}  - E_{7}  +  E_{8},\\
 \partial^{dual}_{2,1} F_3 &=&  E_{3}  - E_{8}  +  E_{9} ,\\
 \partial^{dual}_{2,1} F_4 &=&  E_{4}  - E_{9}  +  E_{10},\\
 \partial^{dual}_{2,1} F_5 &=&  E_{5}  - E_{10} +  E_{11},\\
 \partial^{dual}_{2,1} F_6 &=&  E_{6}  - E_{11} +  E_{12},\\
 \partial^{dual}_{2,1} F_7 &=&  -E_{1}- E_{2}- E_{3}- E_{4}- E_{5}- E_{6}.
\end{eqnarray*}
If we assign coordinate vectors to the dual faces and the dual edges as before
\begin{eqnarray*}
 F_1 &\longleftrightarrow& (1,0,0,0,0,0,0),\\
 F_2 &\longleftrightarrow& (0,1,0,0,0,0,0),\\
 F_3 &\longleftrightarrow& (0,0,1,0,0,0,0),\\
 F_4 &\longleftrightarrow& (0,0,0,1,0,0,0),\\
 F_5 &\longleftrightarrow& (0,0,0,0,1,0,0),\\
 F_6 &\longleftrightarrow& (0,0,0,0,0,1,0),\\
 F_7 &\longleftrightarrow& (0,0,0,0,0,0,1),
\end{eqnarray*}
and 
\begin{eqnarray*}
 E_1 &\longleftrightarrow& (1,0,0,0,0,0,0,0,0,0,0,0),\\
 E_2 &\longleftrightarrow& (0,1,0,0,0,0,0,0,0,0,0,0),\\
 E_3 &\longleftrightarrow& (0,0,1,0,0,0,0,0,0,0,0,0),\\
 E_4 &\longleftrightarrow& (0,0,0,1,0,0,0,0,0,0,0,0),\\
 E_5 &\longleftrightarrow& (0,0,0,0,1,0,0,0,0,0,0,0),\\
 E_6 &\longleftrightarrow& (0,0,0,0,0,1,0,0,0,0,0,0),\\
 E_7 &\longleftrightarrow& (0,0,0,0,0,0,1,0,0,0,0,0),\\
 E_8 &\longleftrightarrow& (0,0,0,0,0,0,0,1,0,0,0,0),\\
 E_9 &\longleftrightarrow& (0,0,0,0,0,0,0,0,1,0,0,0),\\
 E_{10} &\longleftrightarrow& (0,0,0,0,0,0,0,0,0,1,0,0),\\
 E_{11} &\longleftrightarrow& (0,0,0,0,0,0,0,0,0,0,1,0),\\
 E_{12} &\longleftrightarrow& (0,0,0,0,0,0,0,0,0,0,0,1),
\end{eqnarray*}
we have the associated matrix
\[\partial_{2,1}^{dual}=\left(\begin{array}{rrrrrrr}
1 & 0 & 0 & 0 & 0 & 0 & -1\\
0 & 1 & 0 & 0 & 0 & 0 & -1\\
0 & 0 & 1 & 0 & 0 & 0 & -1\\
0 & 0 & 0 & 1 & 0 & 0 & -1\\
0 & 0 & 0 & 0 & 1 & 0 & -1\\
0 & 0 & 0 & 0 & 0 & 1 & -1\\
1 & -1 & 0 & 0 & 0 & 0 & 0\\
0 & 1 & -1 & 0 & 0 & 0 & 0 \\
0 & 0 & 1 & -1 & 0 & 0 & 0 \\
0 & 0 & 0 & 1 & -1 & 0 & 0 \\
0 & 0 & 0 & 0 & 1 & -1& 0 \\
-1 & 0 & 0 & 0 & 0 & 1& 0 
          \end{array}
\right).\]
Notice that  
\[\partial_{2,1}^{dual}=-\left(\partial_{1,0}\right)^T,\]
which arises from the duality between the two meshes \cite[Section 4.5]{Desbrun}.
In general, the discrete differential to be applied is
\[\nabla_{1,2}^{\wedge,dual}=\left(\partial_{2,1}^{dual}\right)^T=-\partial_{1,0}=-\left(\nabla^\wedge_{0,1}\right)^T.\]

\subsection{Discrete Hodge star}
The discretization of the Hodge star $\star$ uses the   
geometrical ideas described in Section \ref{subsec: Hodge star} and the dual mesh.
More precisely, the 2D Hodge star operator rotates a vector $90^\circ$ counterclockwise. 
For the sake of clarity, let us focus on the edge
 $[v_0,v_6]$, its mesh dual $[v_0,v_6]^*$ and its Hodge star image $\star[v_0,v_6]$. They are represented in Figure \ref{fig: edge dual and star}.
\begin{center}
    \includegraphics[width=0.3\textwidth]{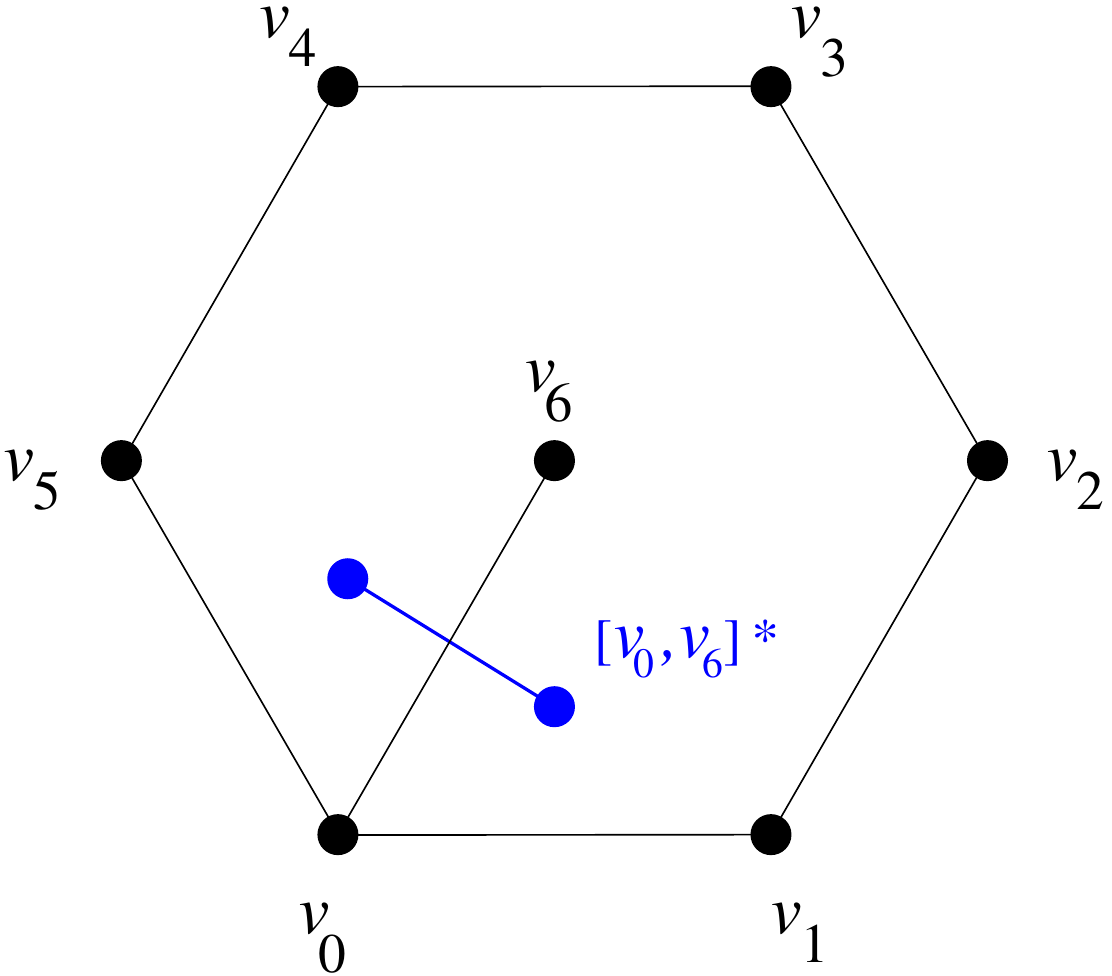}%
    \hspace{.3in}
    \includegraphics[width=0.3\textwidth]{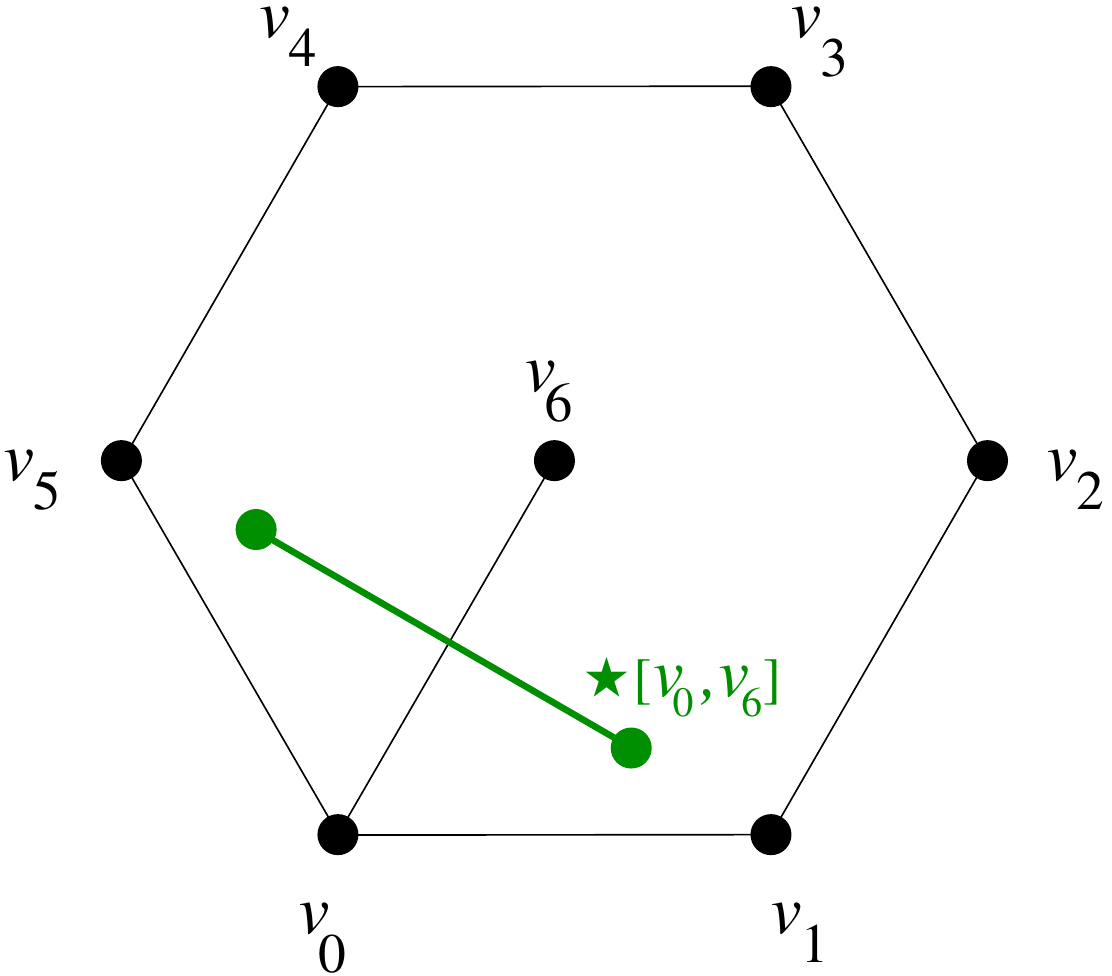}%
    \captionof{figure}{}\label{fig: edge dual and star}
\end{center}
Since
\[{\rm length}(\star[v_0,v_6])={\rm length}([v_0,v_6]),\]
we see that the relationship between the dual edge $[v_0,v_6]^*$ and the geometric $\star[v_0,v_6]$ is the following 
\[{1\over {\rm length}([v_0,v_6]^*)}[v_0,v_6]^*= {1\over {\rm length}([v_0,v_6])}\star[v_0,v_6].\]
As can be seen, when applying the geometric Hodge star to the edges of the mesh, 
we do not end up in the dual mesh but in multiples of the elements of the dual mesh.
Thus, 
$\star[v_0,v_6]$ 
must be scaled to match $[v_0,v_6]^*$.
If we do this to all the Hodge star images, we get
the discrete Hodge star matrix
\[M_{1,1}=\left(\begin{array}{ccccc}
{{\rm length}([v_0,v_6]^*)\over {\rm length}([v_0,v_6])} & 0 & 0 & \dots & 0 \\
0 & {{\rm length}([v_1,v_6]^*)\over {\rm length}([v_1,v_6])} & 0 & \dots & 0 \\
0 & 0 & {{\rm length}([v_2,v_6]^*)\over {\rm length}([v_2,v_6])} & \dots & 0\\
\vdots & \vdots & \vdots &    \ddots &\vdots\\
0 & 0 & 0 & \dots &  {{\rm length}([v_5,v_0]^*)\over {\rm length}([v_5,v_0])}
          \end{array}
\right),\]
where the subindices in $M_{1,1}$ indicate that we are sending 1-dimensional elements of the original mesh to 1-dimensional elements of the dual mesh.

Somewhat less intuitive is the meaning of the geometric Hodge star operator on nodes of the original mesh.
As we saw in Subsection \ref{subsec: Hodge star},
\[\star 1=e_1\wedge e_2,\]
which geometrically means that the Hodge star $\star[v_6]$ must be a polygon with area equal to 1 (classically, it is a parallelogram, but can also be a hexahedron of area 1 as in this example).
    \begin{center}
    \includegraphics[width=0.22\textwidth]{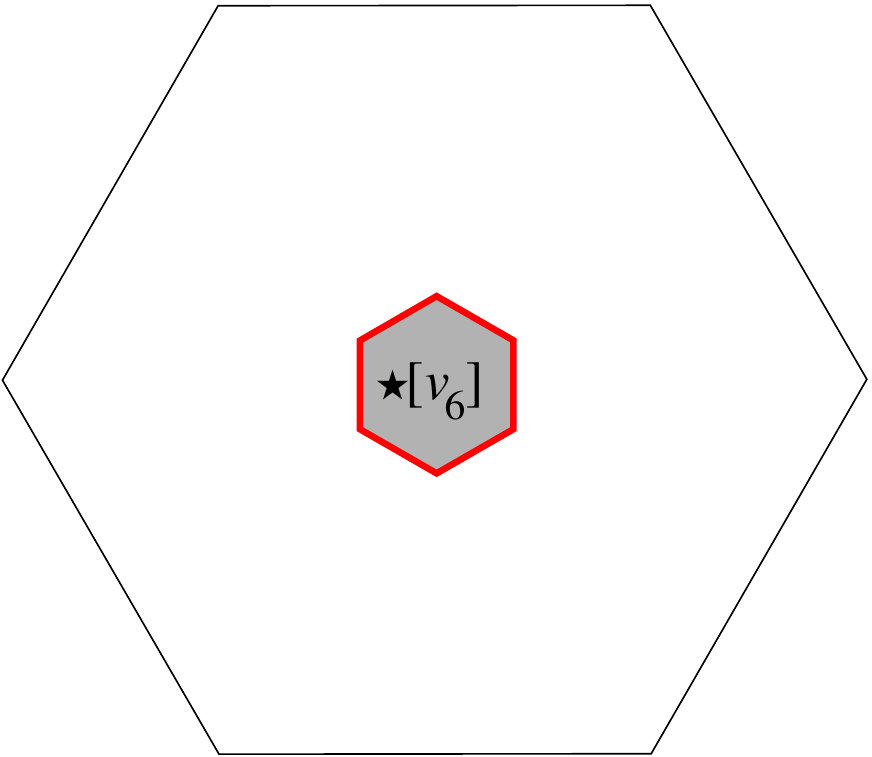}
    \end{center}
However, in the dual mesh we have the polygon $[v_6]^*$ 
    \begin{center}
    \includegraphics[width=0.22\textwidth]{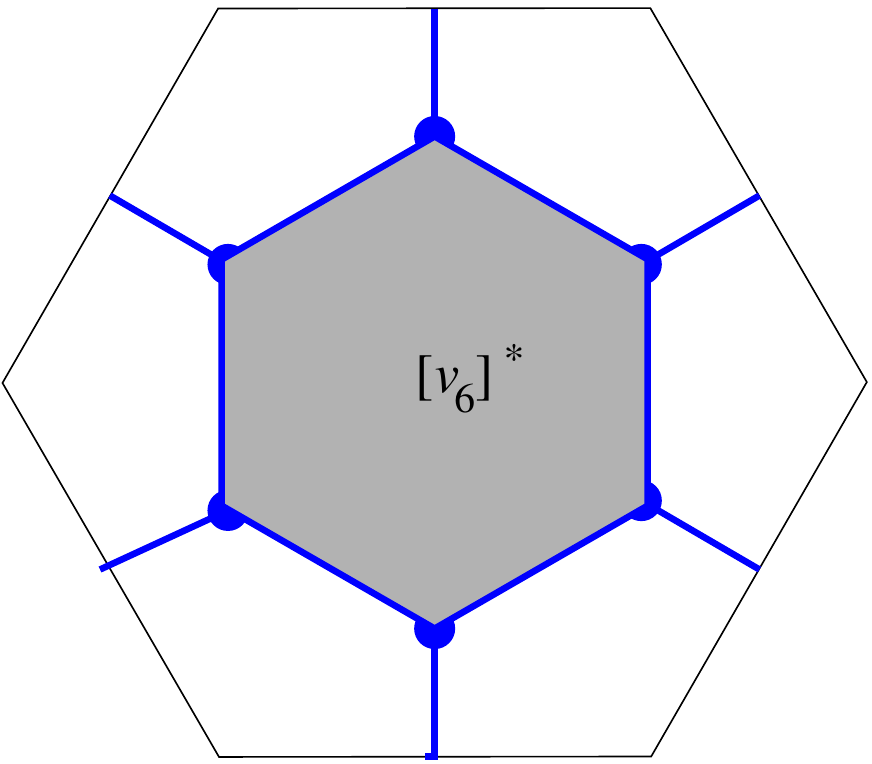}
    \end{center}
so that we need to resize $\star[v_6]$ as follows
\[[v_6]^*={\rm Area}([v_6]^*)\, \star[v_6].\]
If we do that for all the vertices, we obtain the discrete Hodge star matrix 
\[
M_{0,2}=\left(
\begin{array}{cccc}
{\rm Area}([v_0]^*) & 0 & \dots & 0\\
0 & {\rm Area}([v_1]^*) & \dots & 0\\
\vdots & \vdots & \ddots & \vdots\\
0 & 0 & \dots & {\rm Area}([v_6]^*)
\end{array}
\right).\]
The inverse matrix will deal with the case when we take the Hodge star of the 2D polygons in the dual mesh to obtain points (with weight 1) in the original mesh.

In summary, the various matrices $M$ representing the discrete Hodge star operator send elements of the original mesh to elements of the dual mesh.

\subsection{DEC applied to 2D Poisson's equation}

Consider the 2D Poisson's equation 
\[\kappa\Delta f= q.\]
As we have seen, this can be rewritten as
\begin{eqnarray*}
\kappa\,\,\star\nabla^\wedge\star\nabla(f) 
   &=& q.
\end{eqnarray*} 
Suppose that we wish to solve the equation on the meshed domain $K$. 
The equation can them be discretized as the matrix equation
\[\kappa \,\,M_{2,0}\,\, \nabla_{1,2}^{\wedge,dual} \,\, M_{1,1}\,\, \nabla^\wedge_{0,1}\,\, [f]=[q]\]
where $M_{0,2},\nabla_{1,2}^{\wedge,dual}=-\nabla_{0,1}^T,M_{1,1}$ and $\nabla^\wedge_{0,1}$ 
denote the matrices corresponding to the relevant mesh, such as the ones
described in the previous subsections, and
$[f]$ and $[q]$ denote the discretizations of the functions $f$ and $q$ at the nodes/vertices. 
Later on, it will be convenient to work with the  equivalent system 
\[\kappa\,\,(\nabla^\wedge_{0,1})^T\,\,M_{1,1}\,\,\nabla^\wedge_{0,1}[f]= M_{0,2}[q]\]

\section{DEC for general triangulations}\label{sec: general triangle}

Since the boundary operator is really concerned with the connectivity of the mesh and does not change under deformation of the mesh, 
the change in the setup of DEC for a deformed mesh must be contained in the discrete Hodge star matrices. Since 
such matrices are computed in terms of lengths and areas of oriented elements of the mesh, we will now examine how those ingredients 
transform under deformation, a problem that was first considered in \cite{Hirani2}.

\subsection{Dual mesh of an arbitrary triangle}
In order to explain how to implement DEC for general triangulations, let us consider first a mesh consisting of a single well-centered triangle, 
as well as its dual mesh (see
Figure \ref{fig: circumscribed triangle 01}).
    \begin{center}
    \includegraphics[width=0.3\textwidth]{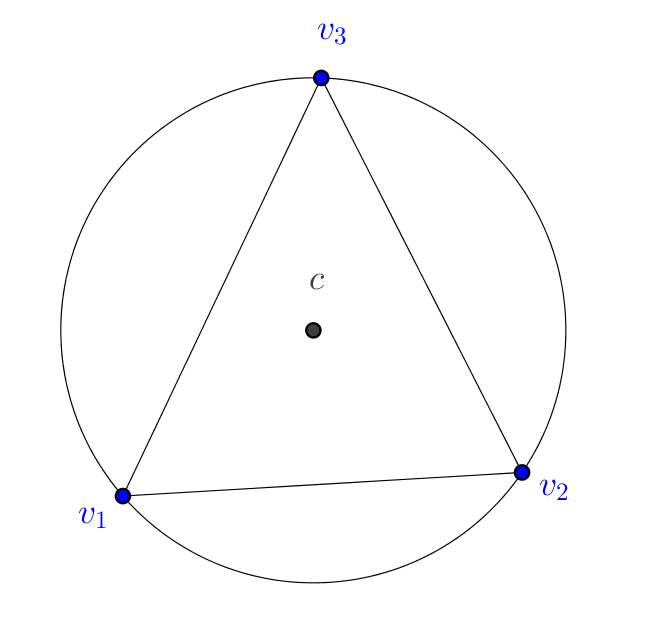}%
    \hspace{.3in}
    \includegraphics[width=0.3\textwidth]{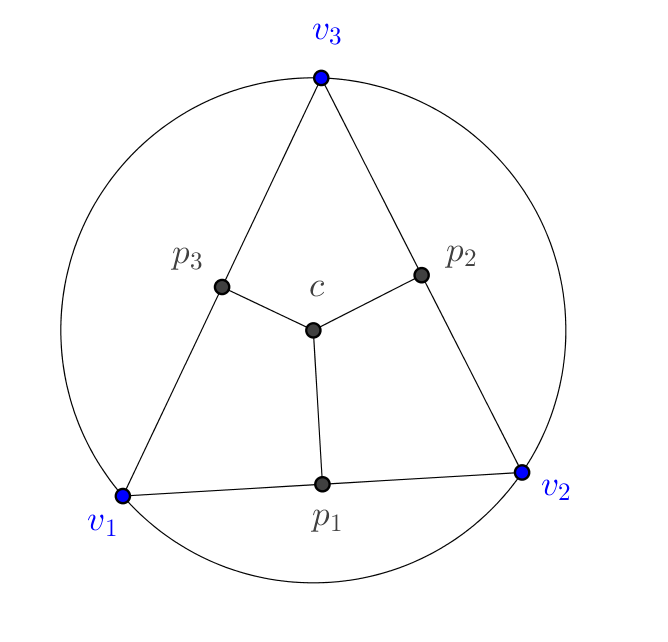}%
    \captionof{figure}{Well-centered triangle and its dual mesh.}\label{fig: circumscribed triangle 01}
    \end{center}
The dual cells are given as follows:
\begin{eqnarray*}
{} [v_1,v_2,v_3]^*&=& [c],\\
{} [v_1,v_2]^*&=& [p_1,c],\\
{} [v_2,v_3]^*&=& [p_2,c],\\
{} [v_3,v_1]^*&=& [p_3,c],\\
{} [v_1]^* &=& [v_1,p_1,c,p_3] ,\\
{} [v_2]^* &=& [v_2,p_2,c,p_1] ,\\
{} [v_3]^* &=& [v_3,p_3,c,p_2] .
\end{eqnarray*}
Now consider the cell $[v_3]^* = [v_3,p_3,c,p_2]$ subdivided as in 
Figure \ref{fig: circumscribed triangle 03}.
    \begin{center}
    \includegraphics[width=0.3\textwidth]{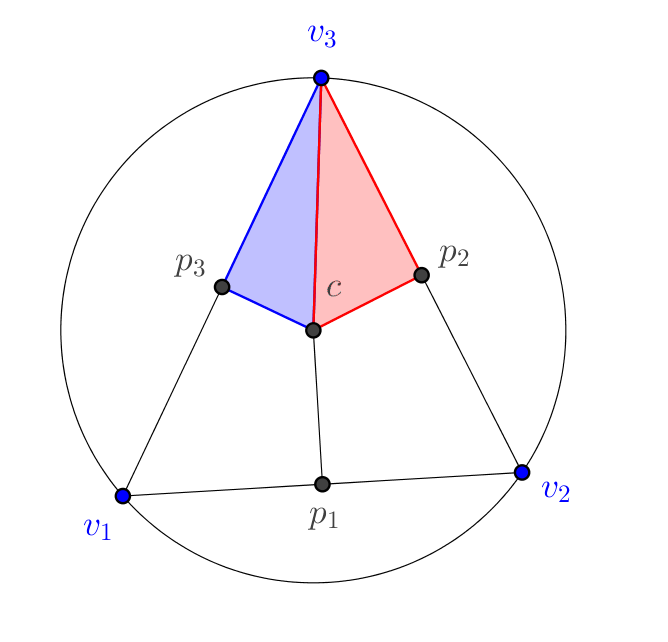}%
    \captionof{figure}{The subdivision of a 2-dimensional dual cell of a well-centered triangle.}\label{fig: circumscribed triangle 03}
    \end{center}
If we deform continuously the triangle $[v_1,v_2,v_3]$ to become an obtuse triangle as in 
Figure \ref{fig: circumscribed triangle 04},
    \begin{center}
    {\small (a)}\includegraphics[width=0.3\textwidth]{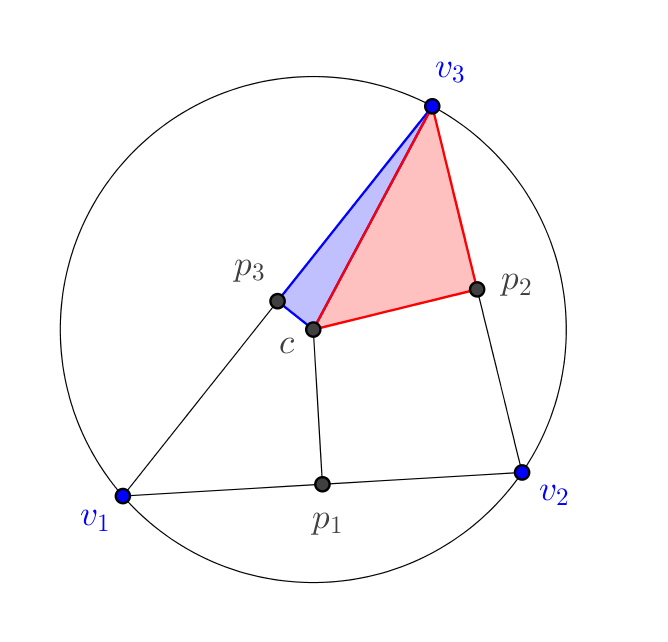}%
    {\small (b)}\includegraphics[width=0.3\textwidth]{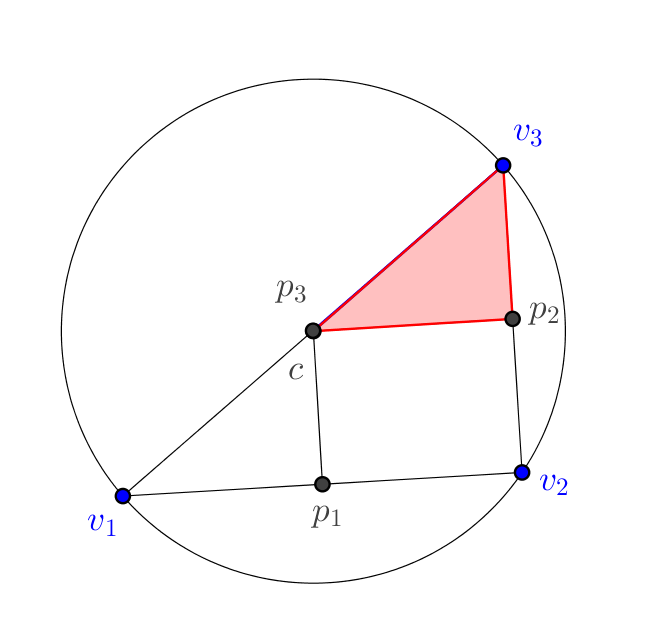}%
    {\small (c)}\includegraphics[width=0.3\textwidth]{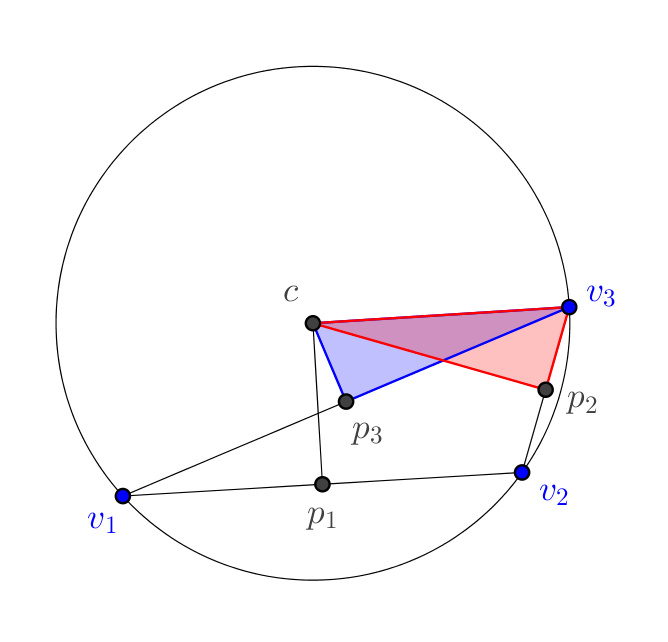}%
    \captionof{figure}{Dual cells of a well-centered triangle.}\label{fig: circumscribed triangle 04}
    \end{center}
we see in Figures \ref{fig: circumscribed triangle 04}(a) and \ref{fig: circumscribed triangle 04}(b) that the area of the blue subtriangle $[v_3,c,p_3]$
decreases to 0 and in Figure \ref{fig: circumscribed triangle 04}(c) that it is completely outside of the triangle and, therefore, must be assigned a negative sign.
The same can be said about the 1-dimensional cell $[p_3,c]$, which originally is completely contained in the triangle $[v_1,v_2,v_3]$, its size reduces to zero as the triangle is deformed (Figures \ref{fig: circumscribed triangle 04}(a) and \ref{fig: circumscribed triangle 04}(b)), 
and eventually it is completely outside the triangle $[v_1,v_2,v_3]$ 
(Figure \ref{fig: circumscribed triangle 04}(c))  and a negative sign must be assigned to it.
On the other hand, part of the red subtriangle $[v_3,c,p_2]$ still intersects the interior of the triangle $[v_1,v_2,v_3]$ and, 
therefore, no assignment of sign is needed. Similarly for the segment $[p_2,c]$.
In terms of numerical simulations, an implementation in terms of determinants contains intrinscally the aforementioned change of signs.

\subsection{Dual mesh of a general triangulation}

Now consider the well-centered mesh and its dual in
Figure \ref{fig: Well-centered hexagonal mesh and its dual mesh}.
    \begin{center}
    {\small (a)}\includegraphics[width=0.3\textwidth]{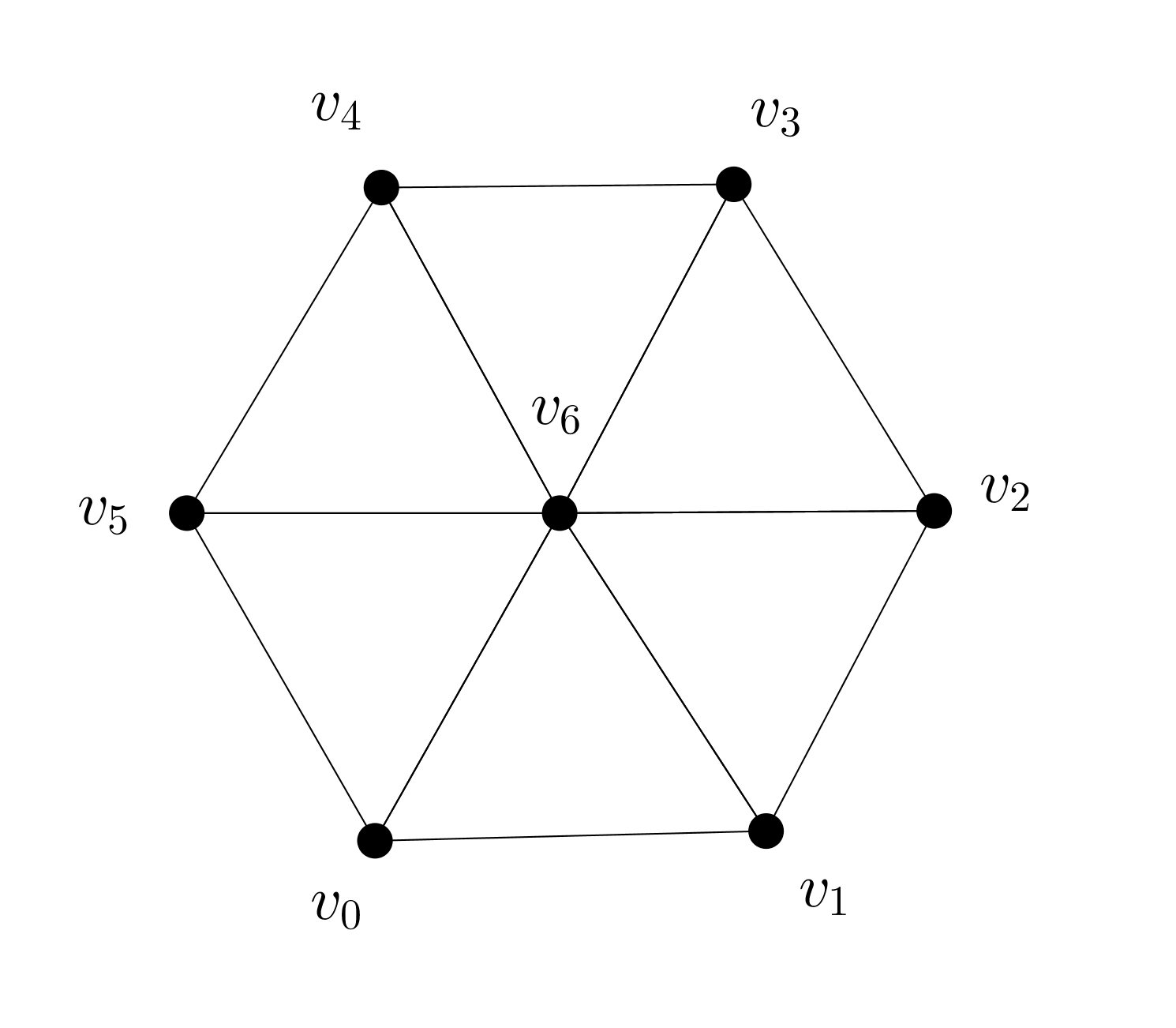}%
    \hspace{.3in}
    {\small (b)}\includegraphics[width=0.29\textwidth]{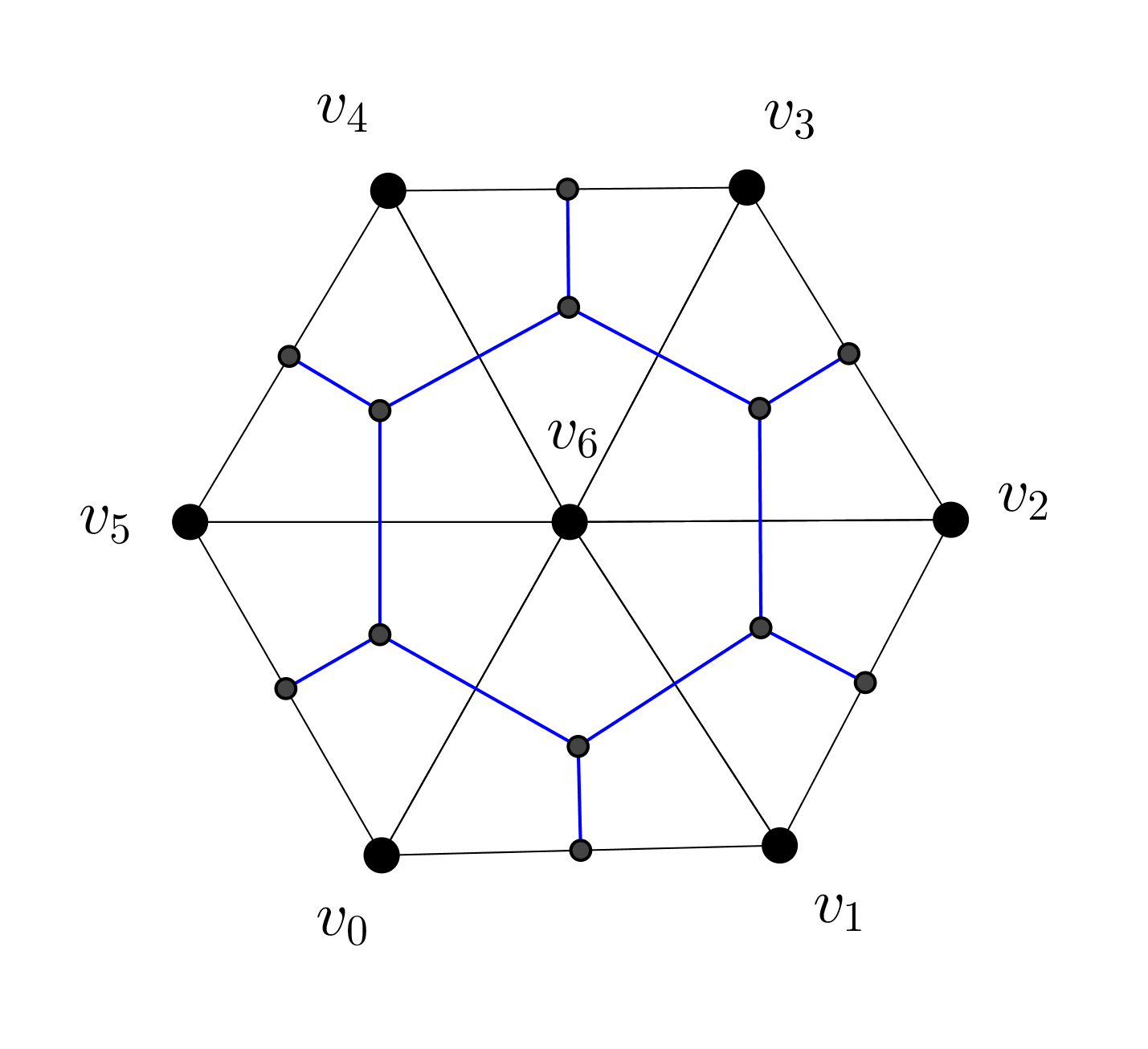}%
    \captionof{figure}{(a) Well-centered triangular mesh of hexagon. (b) Dual mesh.}\label{fig: Well-centered hexagonal mesh and its dual mesh}
    \end{center}
Observe the deformation of the blue-colored dual cell $[v_6]^*$  in
Figure \ref{fig: deformation of dual cell}(a) as the vertex $v_0$ is moved to make the triangle 
$[v_0,v_6,v_5]$ non-well-centered in Figure \ref{fig: deformation of dual cell}(b) and the vertex $v_4$ is moved to make the triangle $[v_4,v_5,v_6]$ 
non-well-centered in Figure \ref{fig: deformation of dual cell}(c). 
    \begin{center}
    {\small (a)}\includegraphics[width=0.3\textwidth]{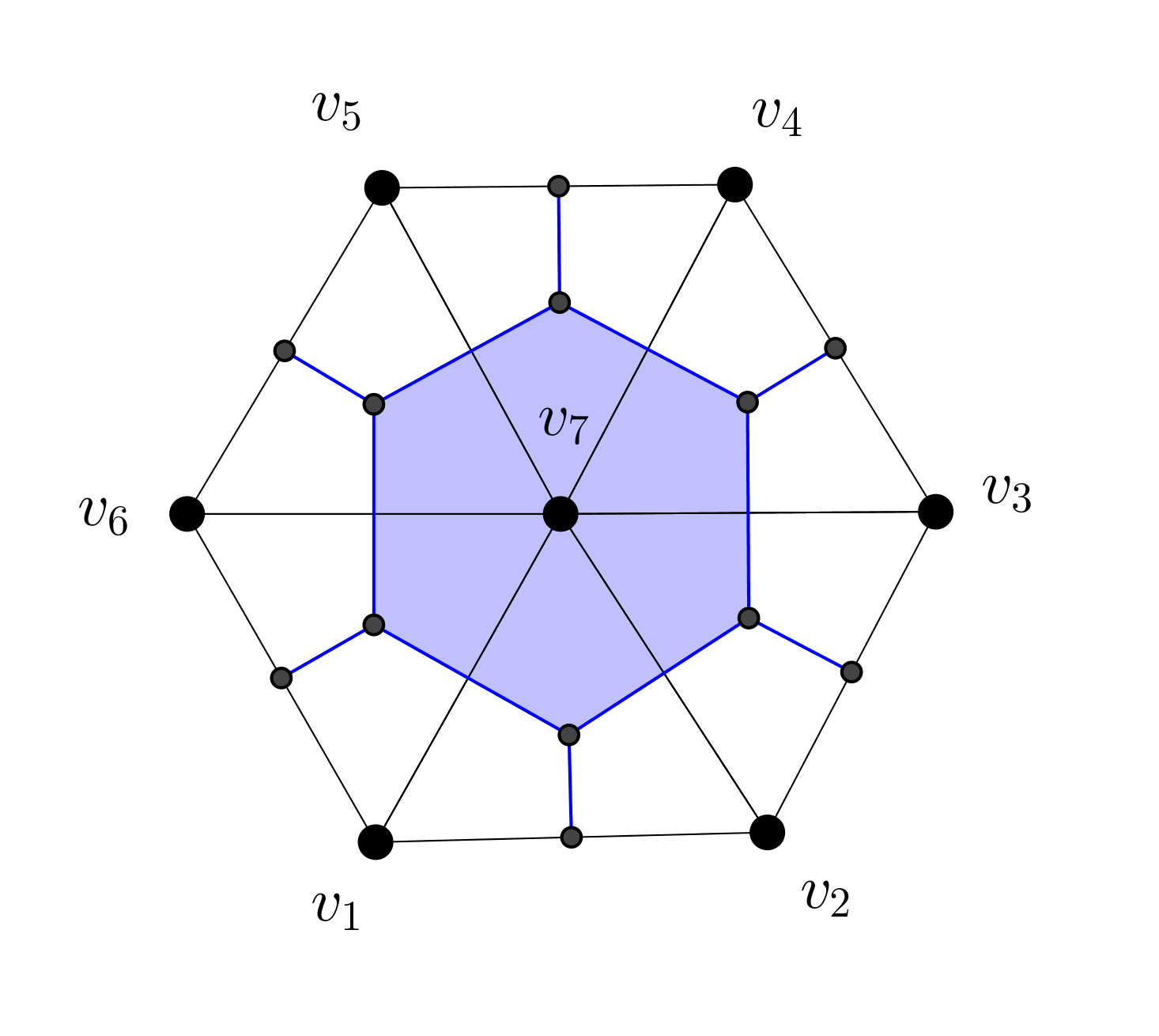}%
    {\small (b)}\includegraphics[width=0.3\textwidth]{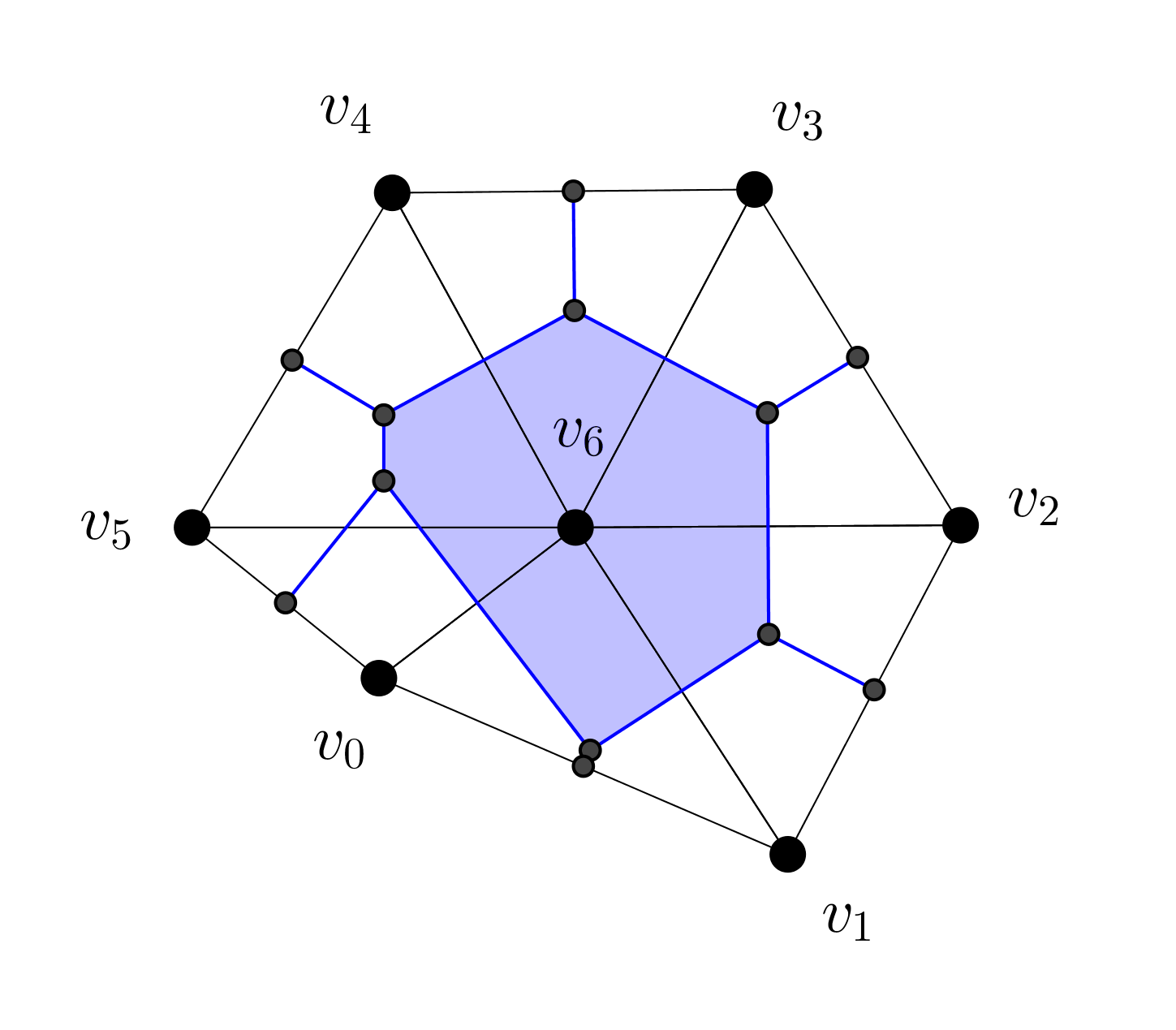}%
    {\small (c)}\includegraphics[width=0.3\textwidth]{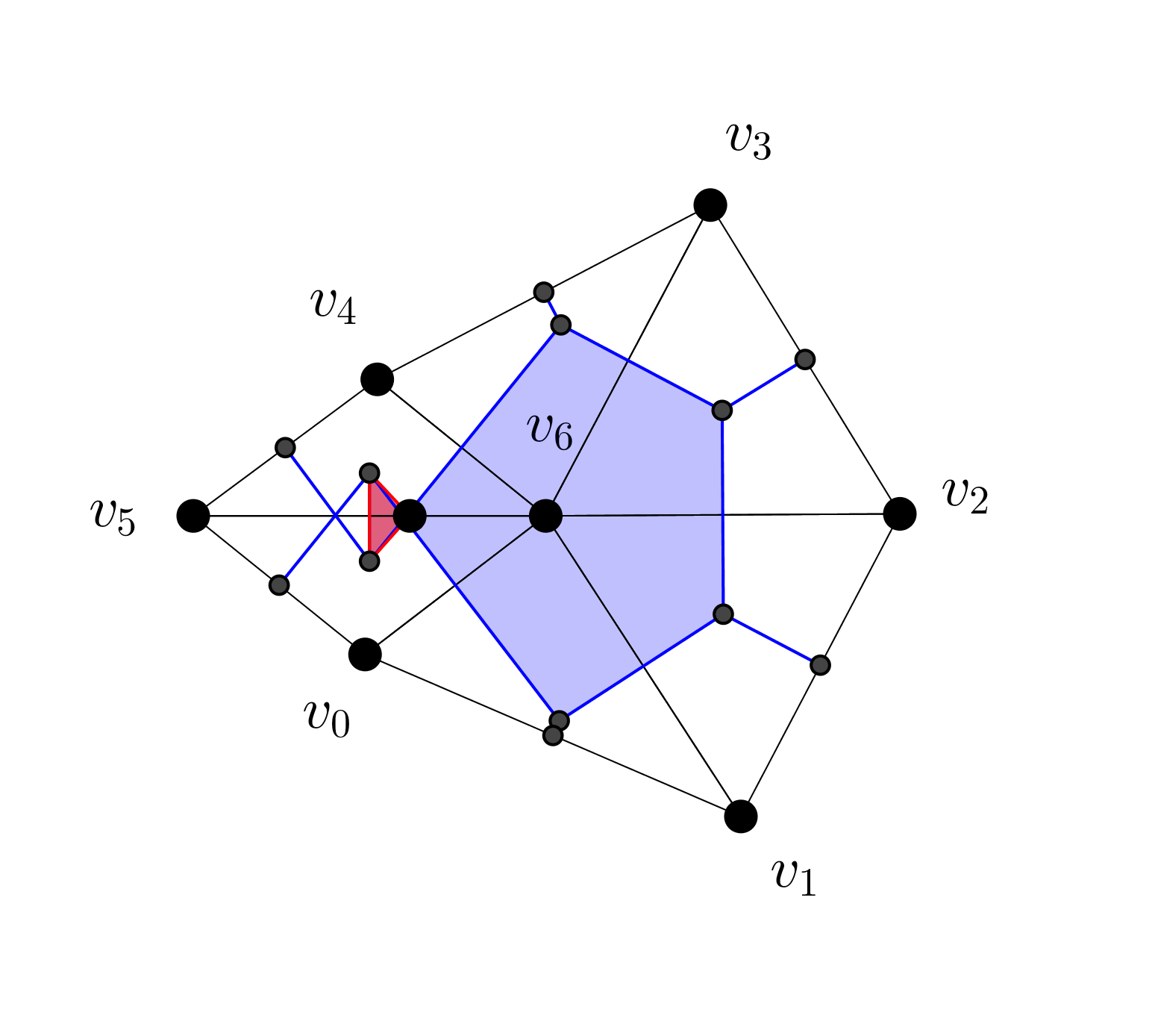}%
    \captionof{figure}{Deformation of a 2-dimensional dual cell as the hexagonal region is deformed.}\label{fig: deformation of dual cell}
    \end{center}
We have colored in red the part of the dual cell $[v_6]^*$ in Figure \ref{fig: deformation of dual cell}(c) that must have the opposite orientation of the blue-colored part.
Similar considerations apply to the edges.

\section{Numerical examples}\label{sec: example}

\subsection{First example}
Let us solve the Poisson equation in a circle of radius one under the following conditions (see Figure \ref{fig: first example geometry}): 
\begin{itemize}
 \item heat difussion constant $k=1$;
 \item source term $q= -1$;
 \item Dirichlet boundary condition $u=10$.
\end{itemize}

The exact solution is
\[
u(x,y)=\frac{1}{4}(1-x^2-y^2)+10.
\]
    \begin{center}
    \includegraphics[width=0.4\textwidth]{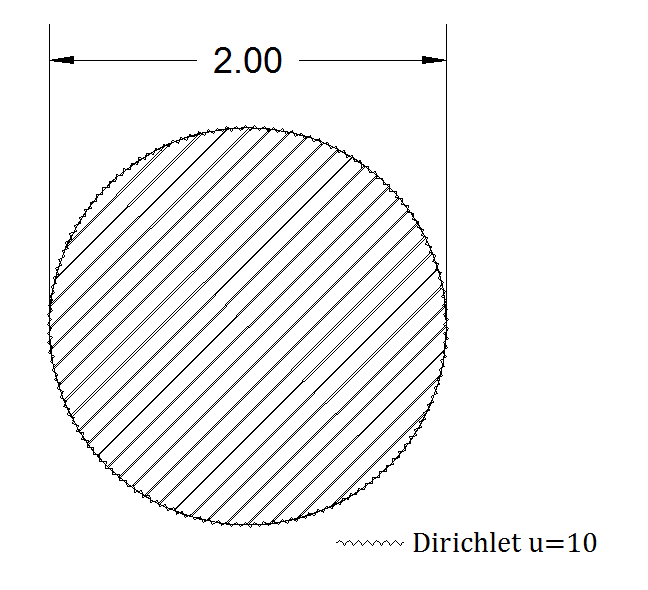}
    \captionof{figure}{Disk of radius one with difussion constant $k=1$, subject to a heat source $q=-1$.}\label{fig: first example geometry}
    \end{center}
The meshes used in this example are shown in Figure \ref{fig: first example meshes} and vary from very coarse to very fine.
    \begin{center}
    {\small (a)}\includegraphics[width=0.25\textwidth]{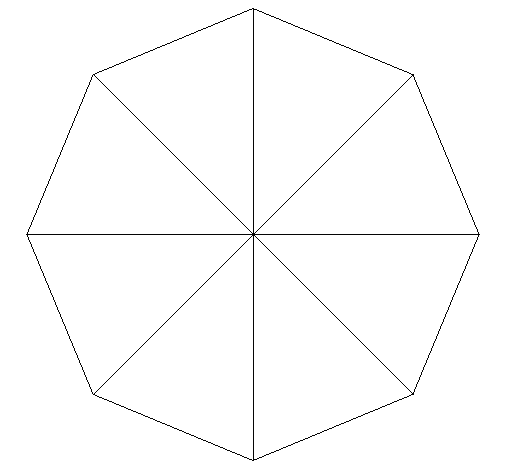}%
    \hspace{.2in}
    {\small (b)}\includegraphics[width=0.25\textwidth]{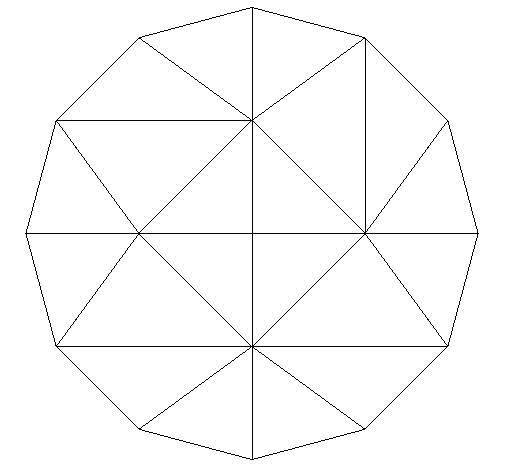}%
    \hspace{.2in}
    {\small (c)}\includegraphics[width=0.25\textwidth]{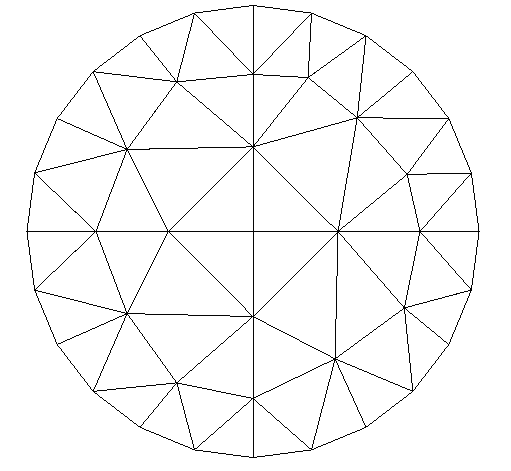}%
    \end{center}
    \begin{center}
    {\small (d)}\includegraphics[width=0.25\textwidth]{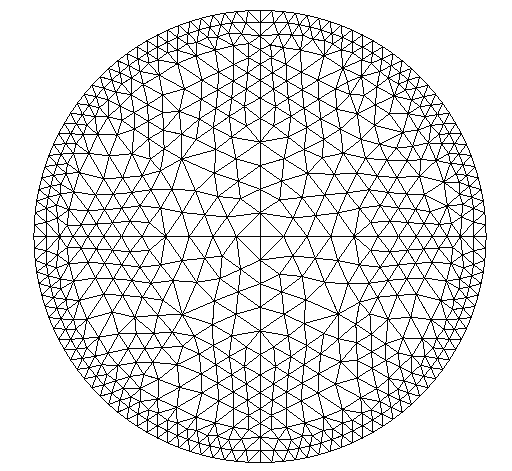}%
    \hspace{.4in}
    {\small (e)}\includegraphics[width=0.25\textwidth]{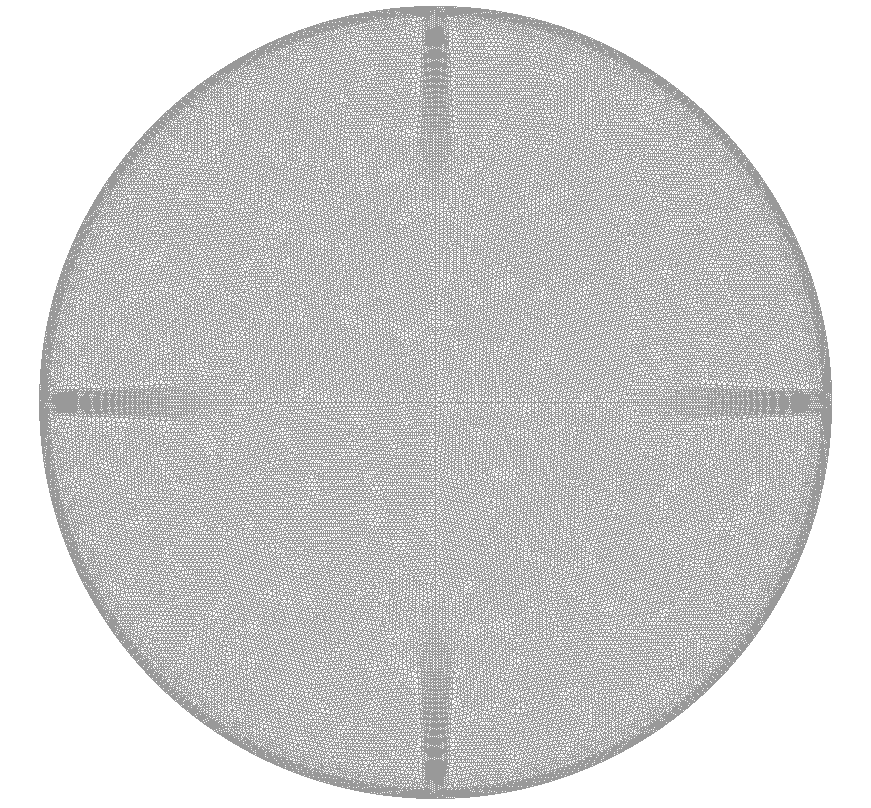}%
    \captionof{figure}{Meshes for unit disk.}\label{fig: first example meshes}
    \end{center}

The numerical results for the maximum temperature value ($u(0,0)=10.25$) are exemplified in Table \ref{table: numerical results} where 
a comparison with the Finite Element Method with linear interpolation functions (FEML) is also shown. 
The FEML methodology that we used in the comparison can be consulted \cite{Zienkiewicz1,Onate,Botello}.
For the sake of completeness in our comparison, here we compute the flux vectors in the same way as in FEML. 
\begin{center}
\[
\begin{array}{|c|c|c|c|c|c|c|}\hline
{\rm Mesh}&	\# {\rm nodes}&	\# {\rm elements}	&	\multicolumn{2}{c|}{\mbox{Max. Temp. Value}}&	\multicolumn{2}{c|}{\mbox{Max. Flux Magnitude}} \\\hline
&&&{\rm DEC} 	&{\rm FEML}&{\rm DEC} 	&{\rm FEML}	 \\\hline
 {\rm Figure \,\,\,\ref{fig: first example meshes}(a)}&	9&	8&	10.250&				10.285& 0.270 & 0.307\\
{\rm Figure \,\,\,\ref{fig: first example meshes}(b)}&	17&	20&	10.250&			10.237& 0.388 & 0.405\\
{\rm Figure \,\,\,\ref{fig: first example meshes}(c)}&	41&	56&	10.250&				10.246& 0.449 & 0.453\\
{\rm Figure \,\,\,\ref{fig: first example meshes}(d)}&	713&	1304&	10.250&				10.250& 0.491 & 0.492\\
{\rm Figure \,\,\,\ref{fig: first example meshes}(e)}&	42298&	83346&	10.250&				10.250& 0.496 & 0.496\\\hline
\end{array}    
\]    
\captionof{table}{Maximum temperature and Flux magnitude values in the numerical simulation.}\label{table: numerical results} 
\end{center}
The temperature distribution and Flux magnitude fields for the finest mesh are shown in Figure \ref{fig: first example temperature field}.
    \begin{center}
    \includegraphics[width=0.35\textwidth]{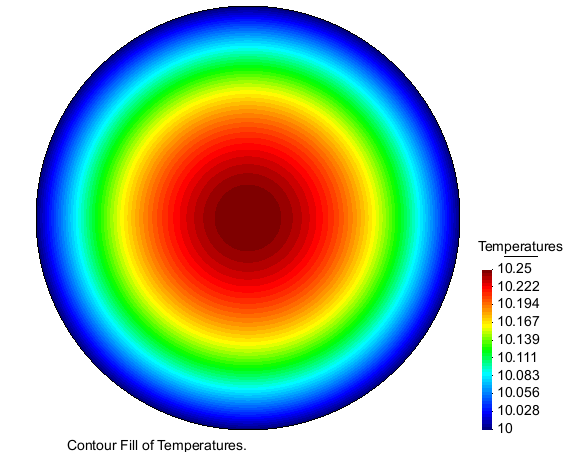}%
    \hspace{.2in}
    \includegraphics[width=0.41\textwidth]{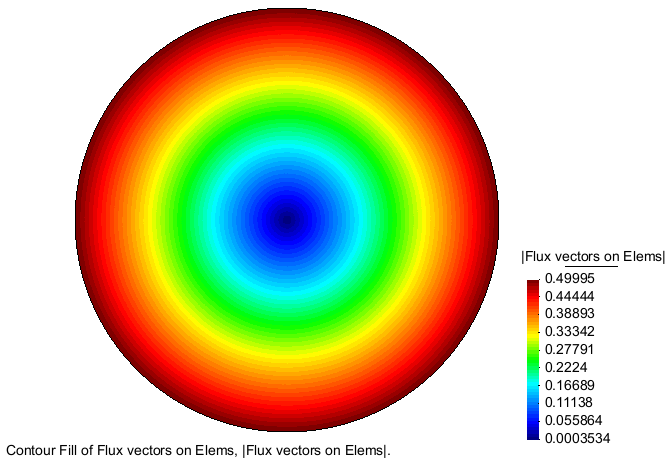}%
    \captionof{figure}{Temperature distribution and Flux magnitude fields for the finest mesh.}\label{fig: first example temperature field}
    \end{center}
Figures \ref{fig: first example diametral graphs}(a), \ref{fig: first example diametral graphs}(b) and \ref{fig: first example diametral graphs}(c) 
show the graphs of the temperature and flux magnitude values along a diameter of the circle for the different meshes of Figures \ref{fig: first example meshes}(b), \ref{fig: first example meshes}(c)
and \ref{fig: first example meshes}(d) respectively.
    \begin{center}
    {\small (a)}\includegraphics[width=0.43\textwidth]{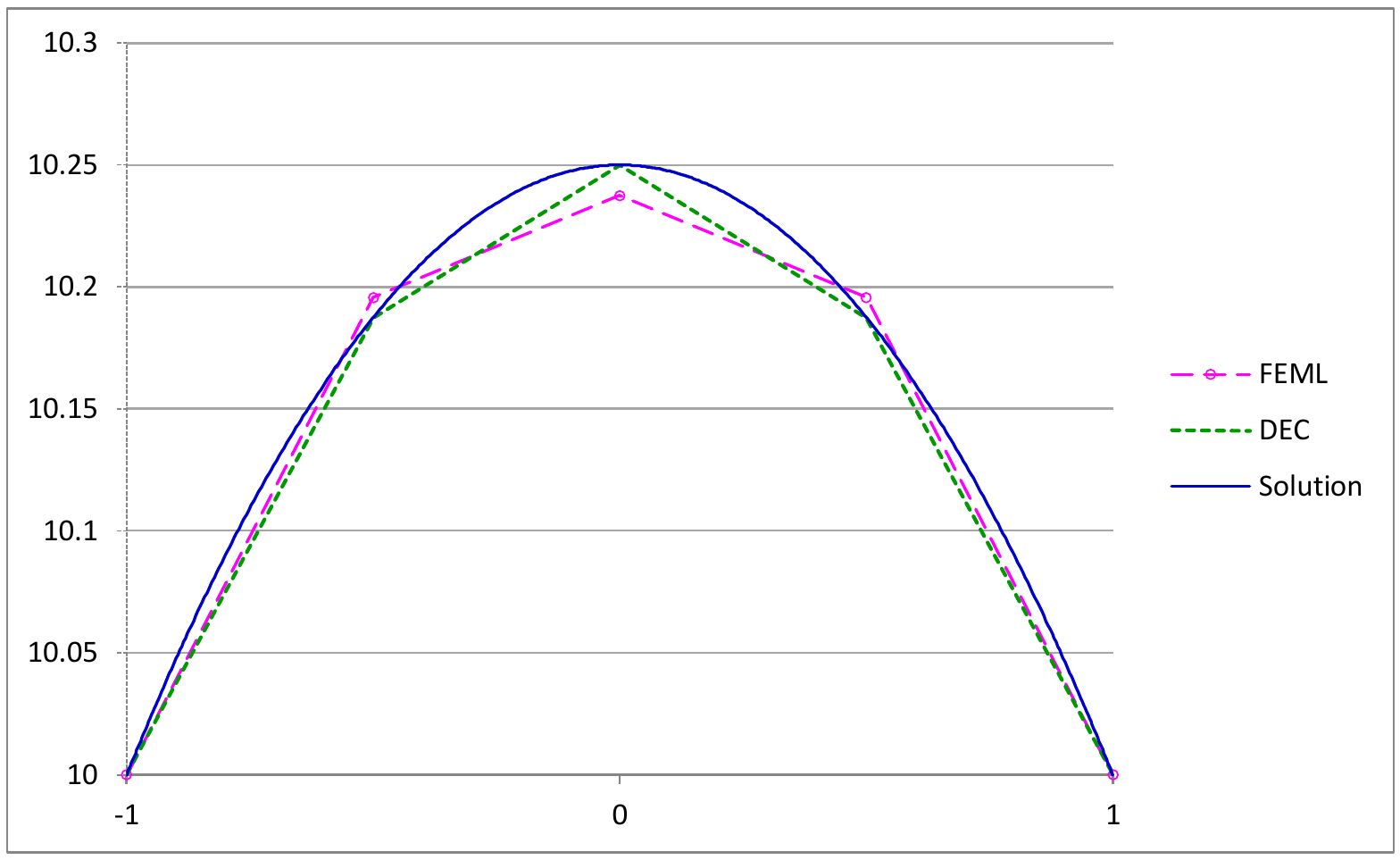}%
    \hspace{.1in}
    \includegraphics[width=0.43\textwidth]{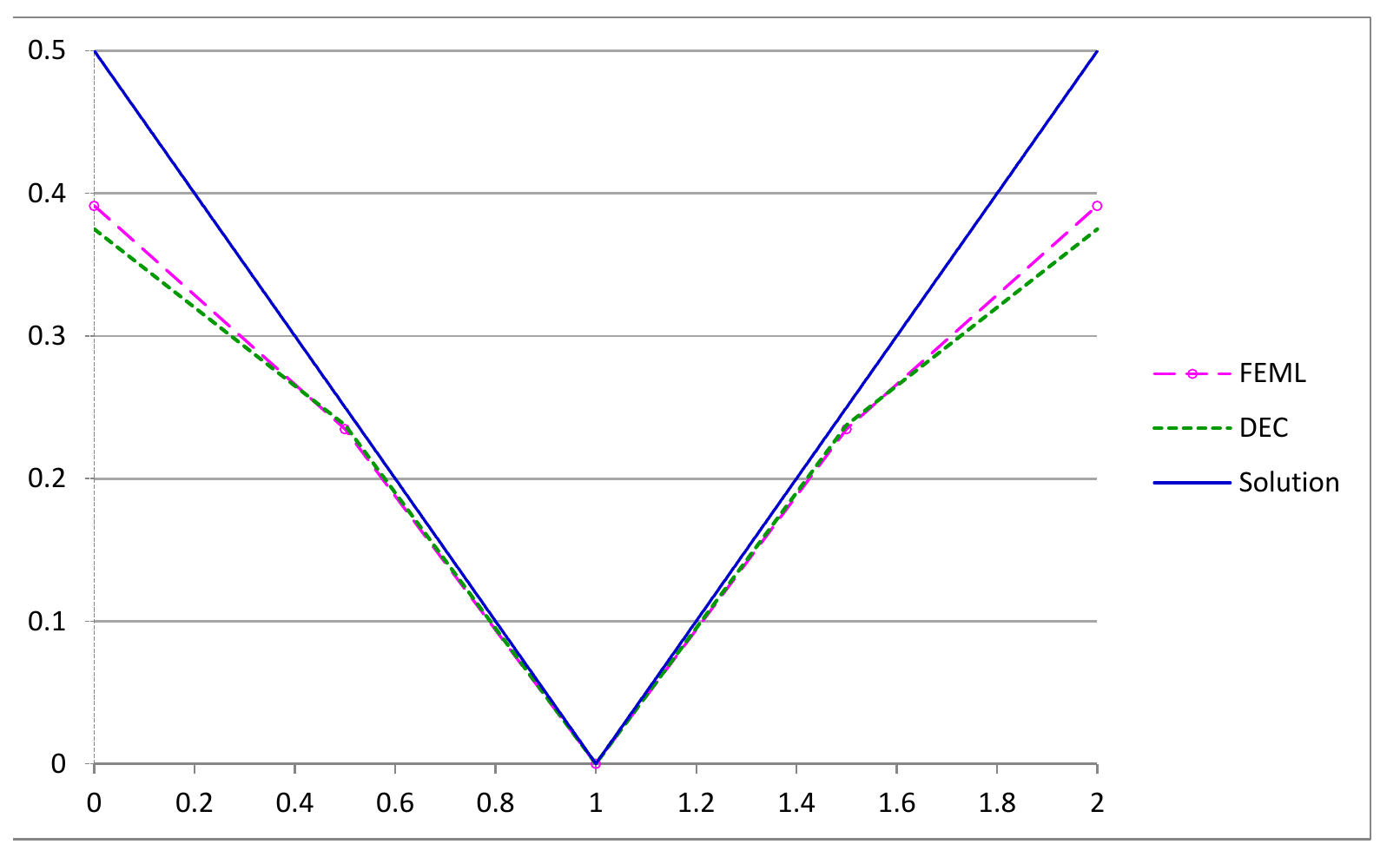}
    \newline
    {\small (b)}\includegraphics[width=0.43\textwidth]{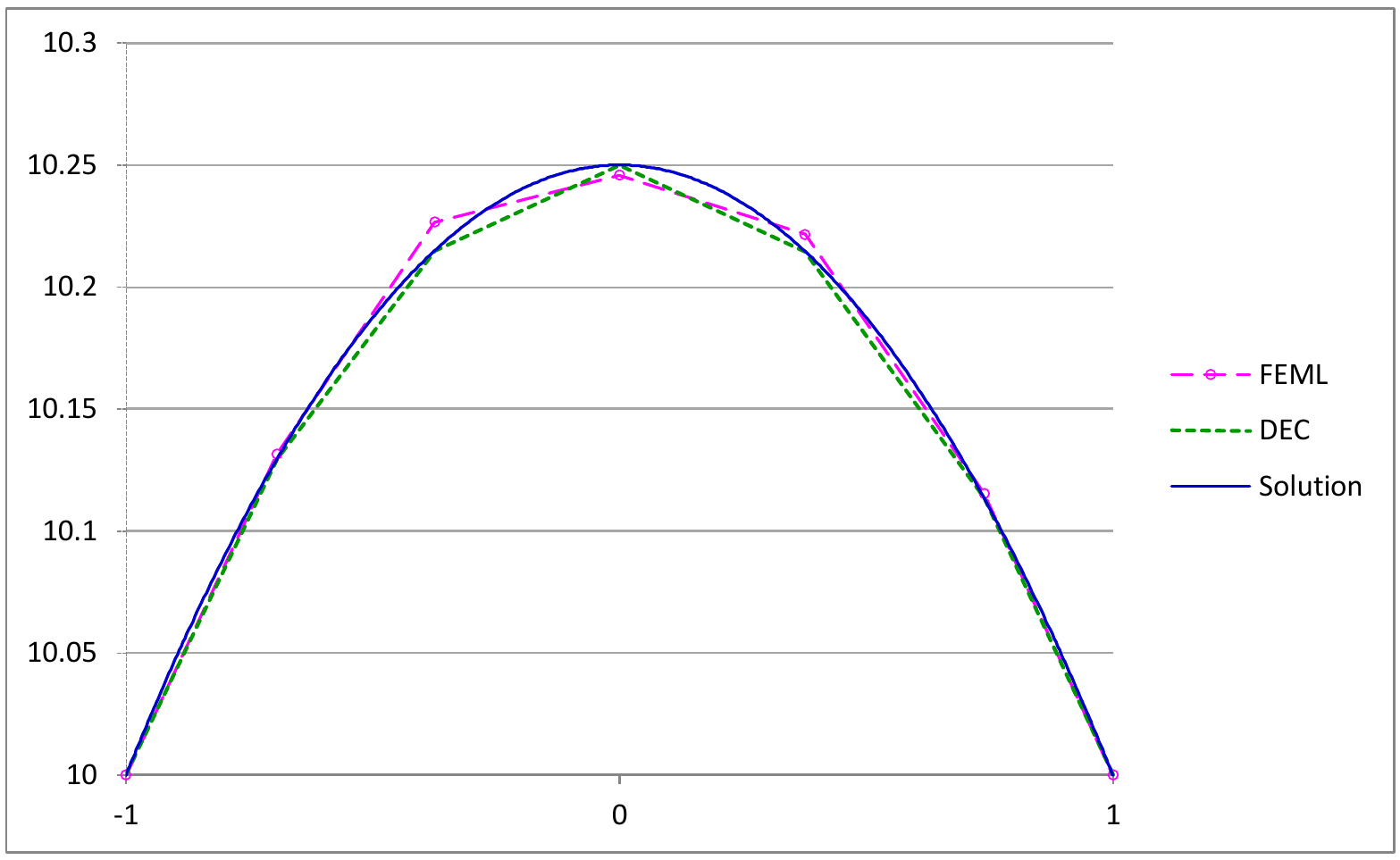}%
    \hspace{.1in}
    \includegraphics[width=0.43\textwidth]{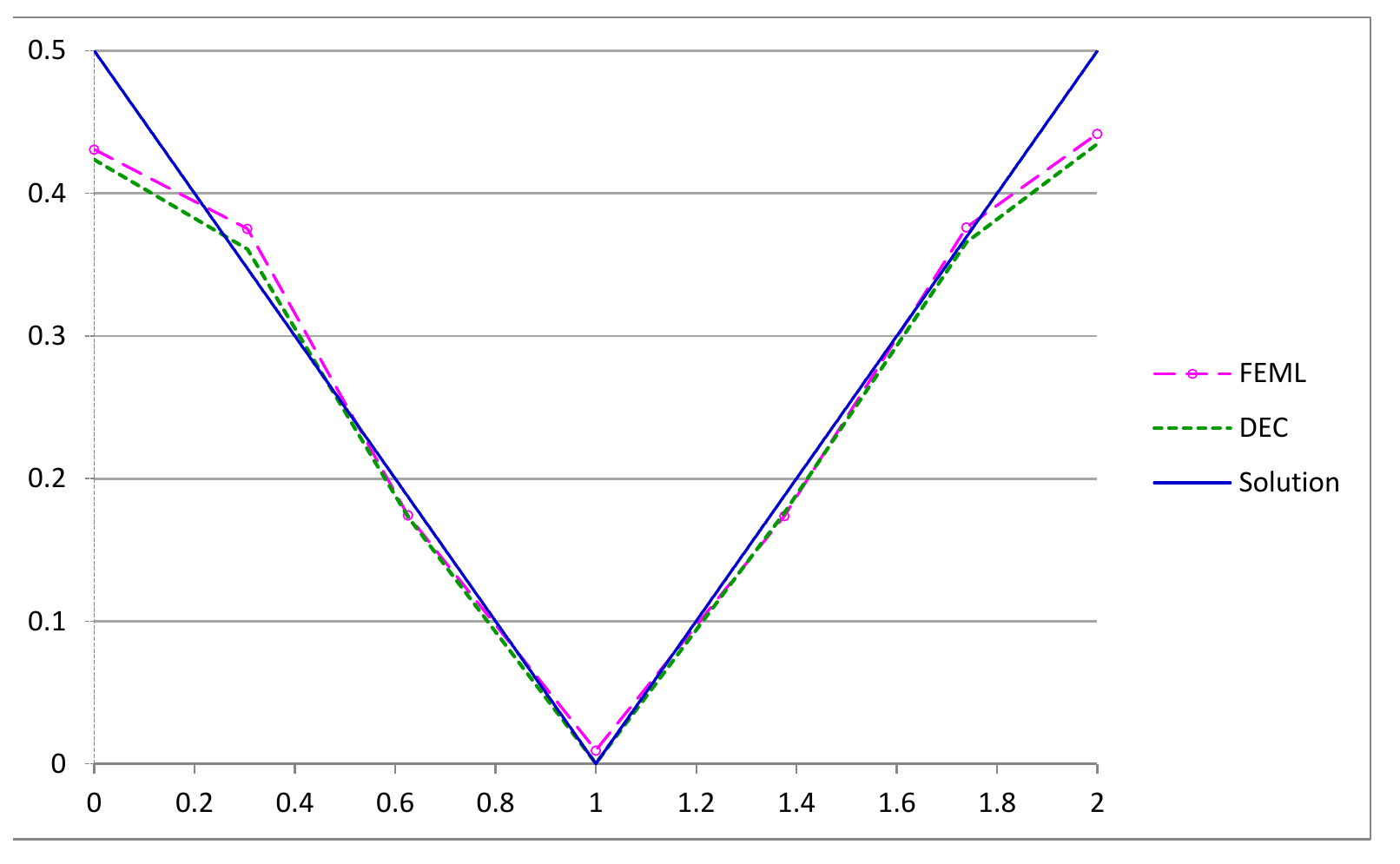}%
    \newline
    \noindent{\small (c)}\includegraphics[width=0.43\textwidth]{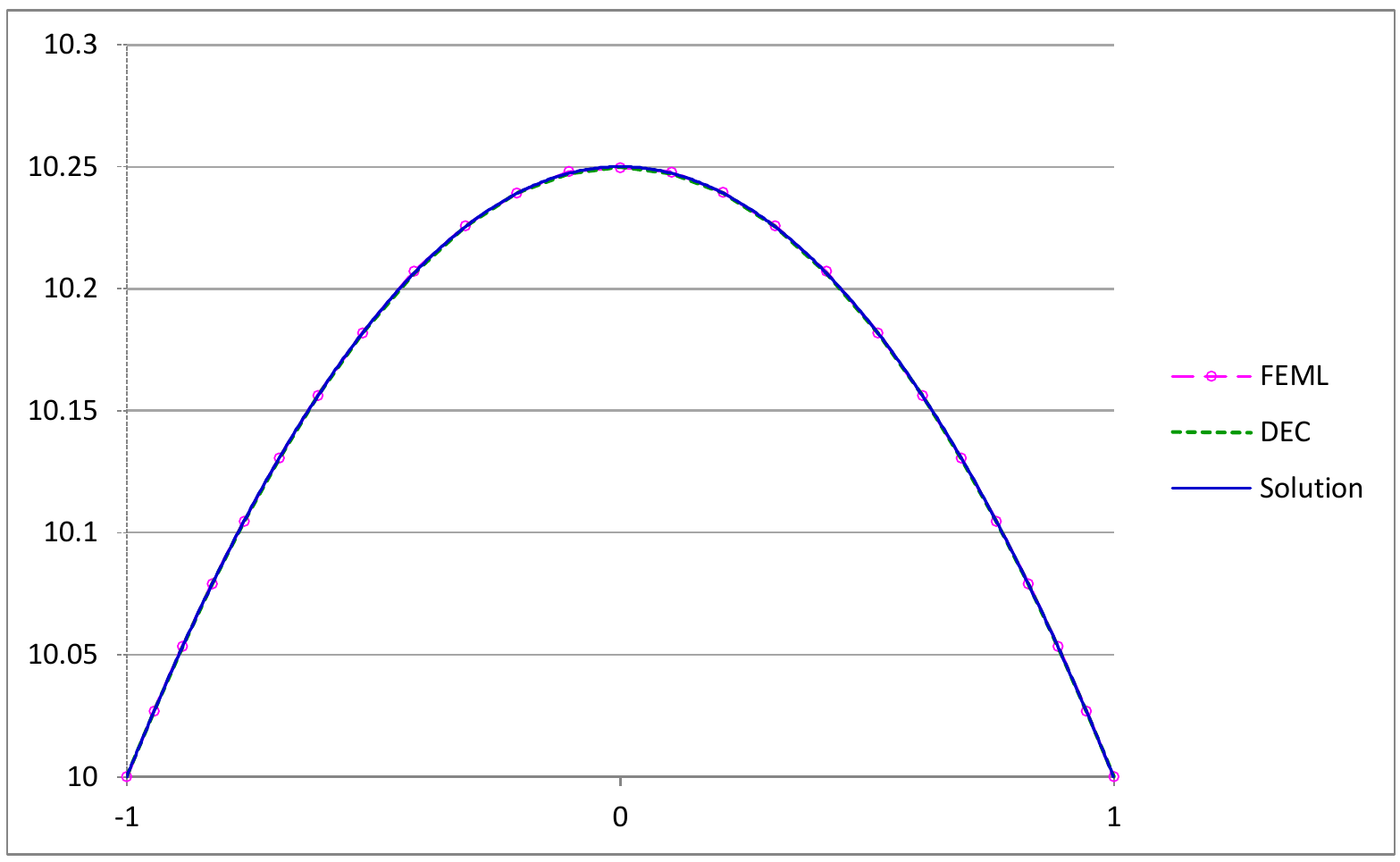}%
    \hspace{.1in}
    \includegraphics[width=0.43\textwidth]{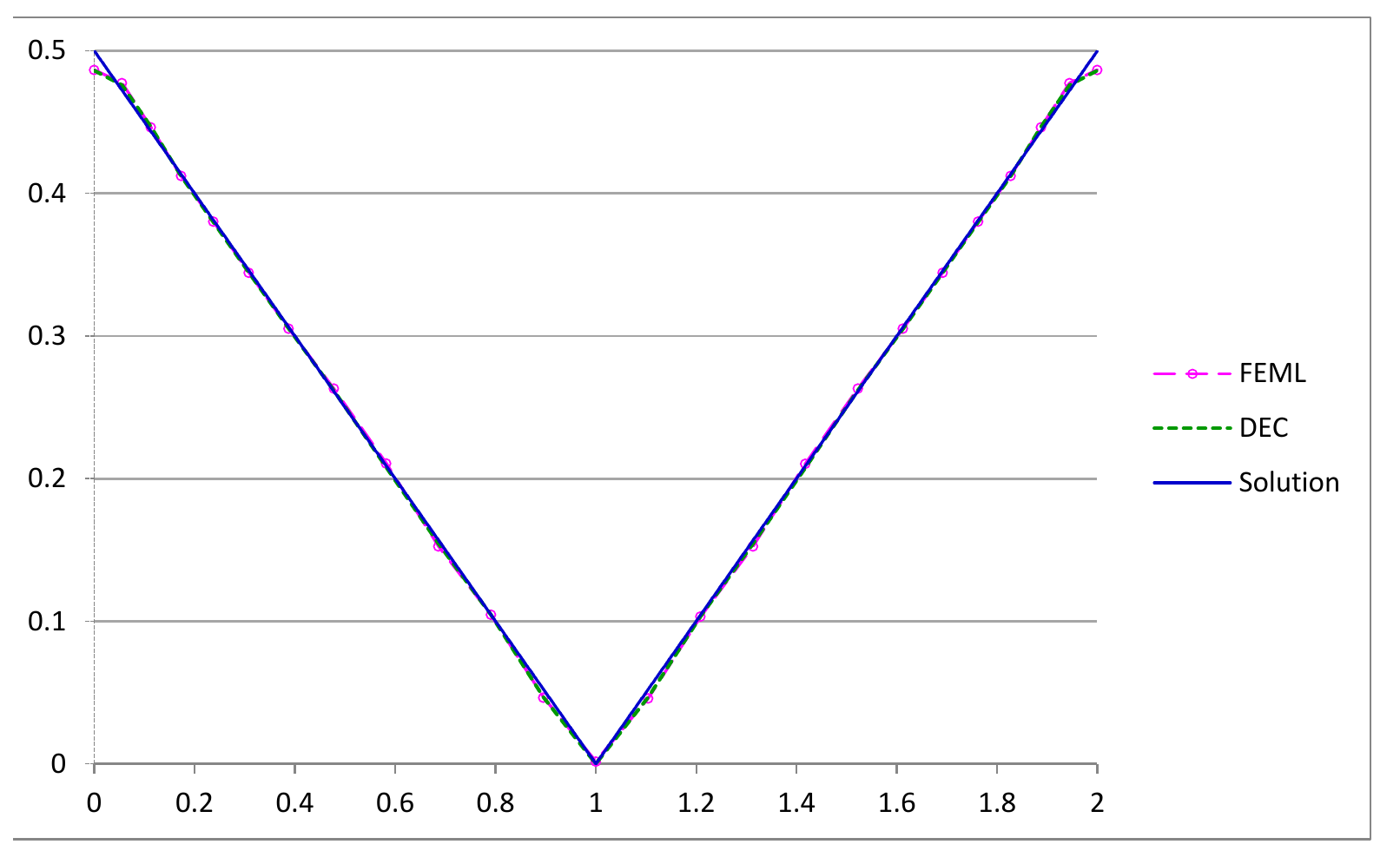}%
    \captionof{figure}{Temperature and Flux magnitude graphs along a diameter of the circle for different meshes: 
    (a) Graphs for the Mesh in Figure \ref{fig: first example meshes}(b); 
    (b) Graphs for the Mesh in Figure \ref{fig: first example meshes}(c); 
    (c) Graphs for the Mesh in Figure \ref{fig: first example meshes}(d); 
    }\label{fig: first example diametral graphs}
    \end{center}

As can be seen from Table \ref{table: numerical results} and Figure \ref{fig: first example diametral graphs}, DEC behaves very well
on coarse meshes. Note that the maximum temperature values calculated with DEC matches the exact 
theoretical value even on coarse meshes.
As expected, the results of DEC and FEML are similar for fine meshes. 
We would also like to point out the the computational costs of DEC and FEML are very similar.

\subsection{Second example}

In this example, we consider a region in the plane whose boundary consists of segments of a straight line, a circle, a parabola, a cubic and an ellipse (see Figure \ref{fig: second example geometry}).
    \begin{center}
    \includegraphics[width=0.39\textwidth]{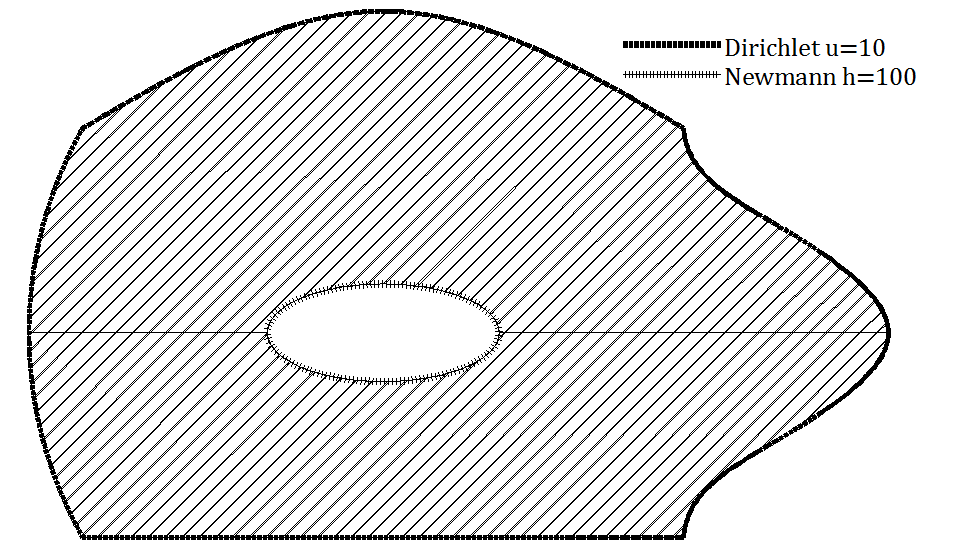}%
     \hspace{.3in}
    \includegraphics[width=0.39\textwidth]{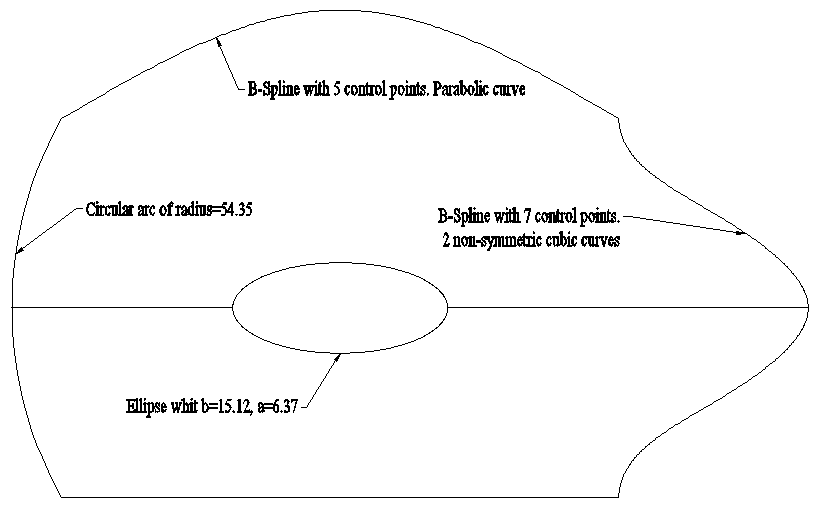}%
    \captionof{figure}{Region with linear, quadratic and cubic boundary segments, together with boundary conditions. 
    }\label{fig: second example geometry}
    \end{center}

The meshes used in this example are shown in Figure \ref{fig: second example meshes} and vary from coarse to very fine.
    \begin{center}
    (a)\includegraphics[width=0.25\textwidth]{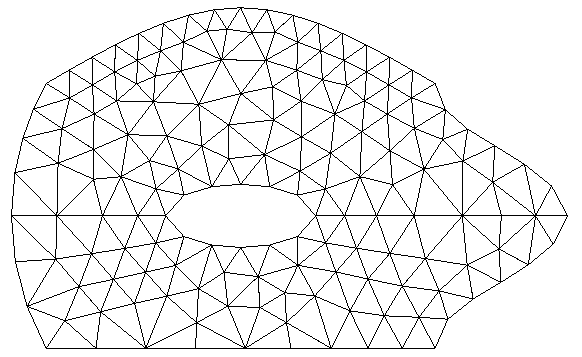}%
    \hspace{.2in}
    (b)\includegraphics[width=0.25\textwidth]{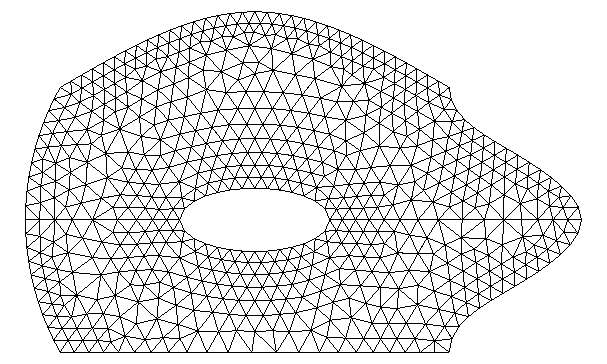}%
    \hspace{.2in}
    (c)\includegraphics[width=0.25\textwidth]{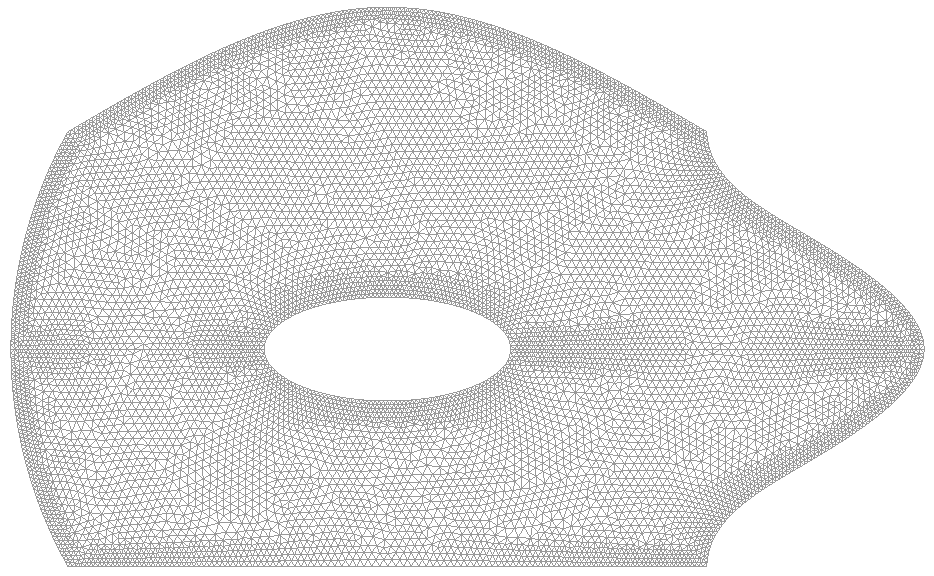}%
    \captionof{figure}{Three of the meshes used in the first example. 
    }\label{fig: second example meshes}
    \end{center}
We set
\begin{itemize}
 \item the source term $q=20.2$; 
 \item the heat diffusion constant $k=80.2$;
 \item Newmann boundary condition $h=100$ along the inner elliptical boundary; 
 \item Dirichlet boundary condition $u=10$ along the external boundary. 
\end{itemize}
The numerical results for the maximum temperature and flux magnitude values are exemplified in Table \ref{table: third example numerical results}.
\begin{center}
\[
\begin{array}{|c|c|c|c|c|c|c|}\hline
&	 &	&	\multicolumn{2}{c|}{\mbox{Max. Temp. Value}}&	\multicolumn{2}{c|}{\mbox{Max. Flux Magnitude}} \\\hline
{\rm Mesh}	& \# {\rm  nodes} &\# {\rm elements}	&{\rm DEC}		&{\rm FEML} & {\rm DEC}  & {\rm FEML}	\\\hline
\mbox{Figure \ref{fig: second example meshes}(a)}	&162	&268		&129.07		&128.92	& 467.22  & 467.25\\
\mbox{Figure \ref{fig: second example meshes}(b)}	&678	&1223		&129.71		&129.70	& 547.87  & 548.42 \\
\mbox{Figure \ref{fig: second example meshes}(c)}	&9489	&18284		&130.00		&130.00	& 591.58  & 591.63  	\\
&27453	&53532		&130.01		&130.01	& 597.37  & 597.38 	\\
&651406	&1295960	&130.02		&130.02	& 602.19  & 602.19 	\\\hline
\end{array}
\]
\captionof{table}{Numerical simulation results.}\label{table: third example numerical results}  
\end{center}
The temperature and flux-magnitude distribution fields are shown in Figure \ref{fig: third example temperature field}.    
    \begin{center}
    (a) \includegraphics[width=0.4\textwidth]{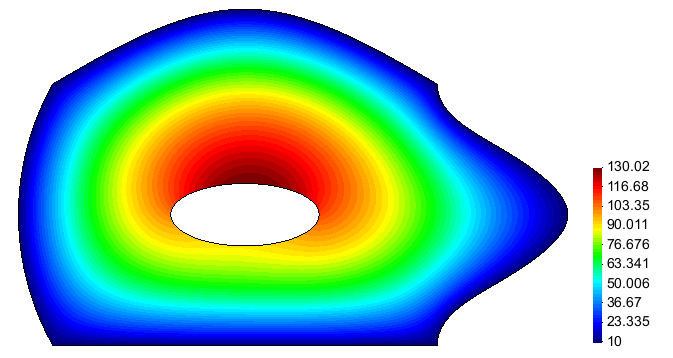}%
    \hspace{.2in}
    (b) \includegraphics[width=0.4\textwidth]{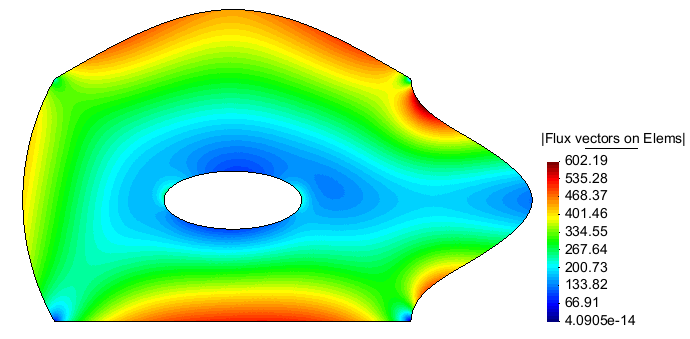}%
    \captionof{figure}{Temperature and Flux magnitude distribution fields.}\label{fig: third example temperature field}
    \end{center}
Figure \ref{fig: third example diametral graphs} shows the graphs of the temperature and the flux magnitude along a horizontal line crossing the elliptical boundary for the first two meshes.
    \begin{center}
    (a)\includegraphics[width=0.4\textwidth]{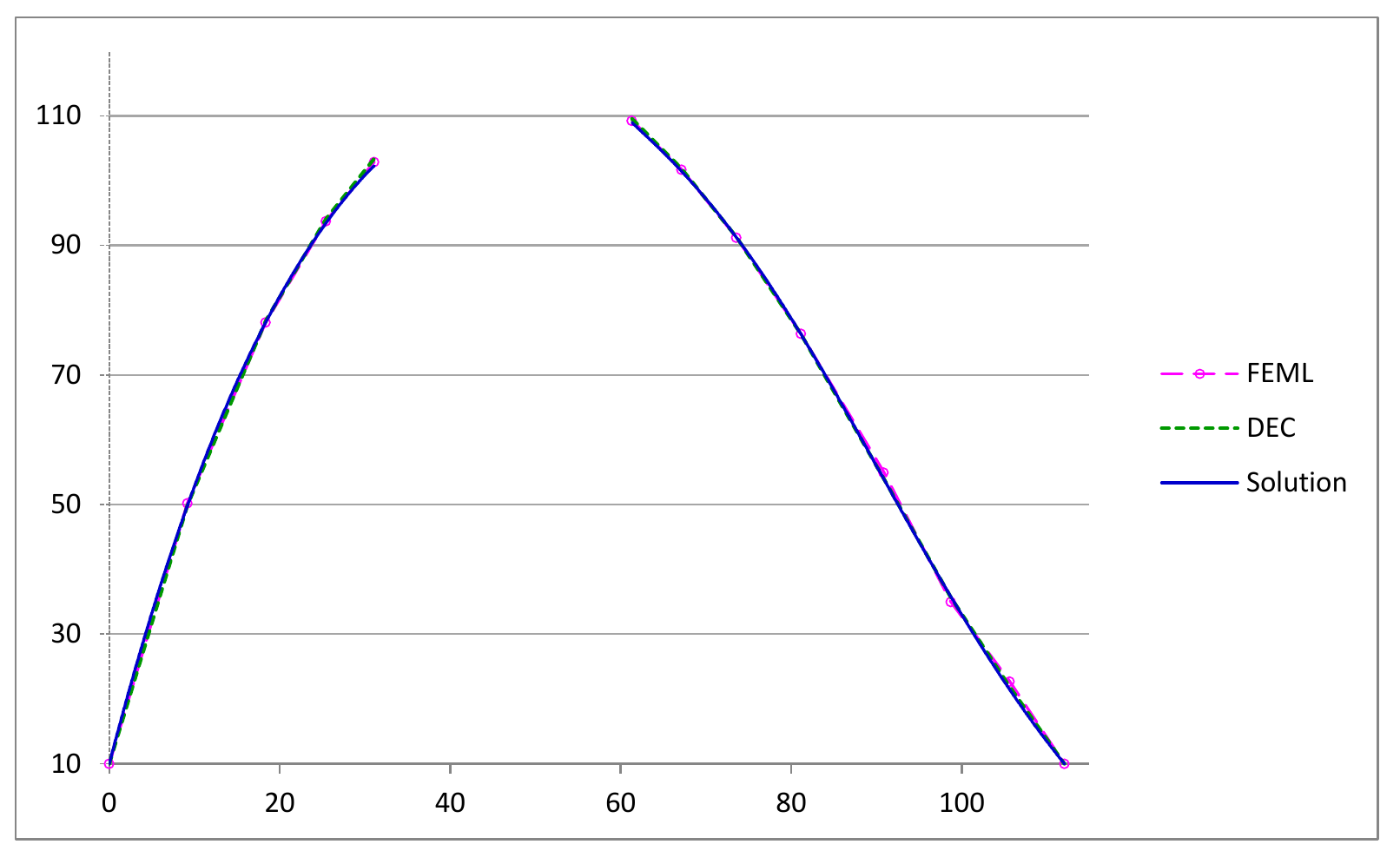}%
    \hspace{.1in}
    \includegraphics[width=0.4\textwidth]{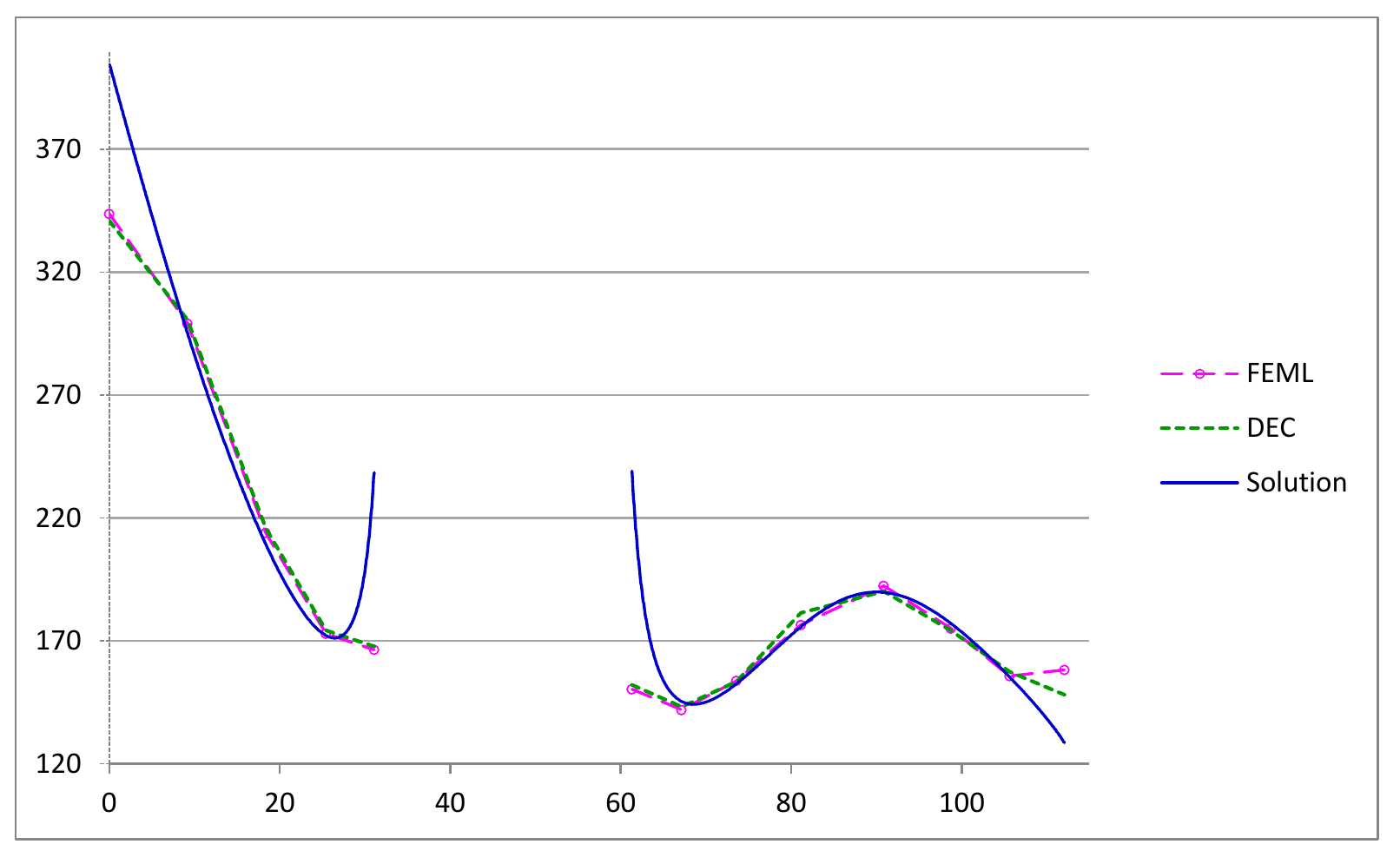}%
    \newline
    \noindent(b)\includegraphics[width=0.4\textwidth]{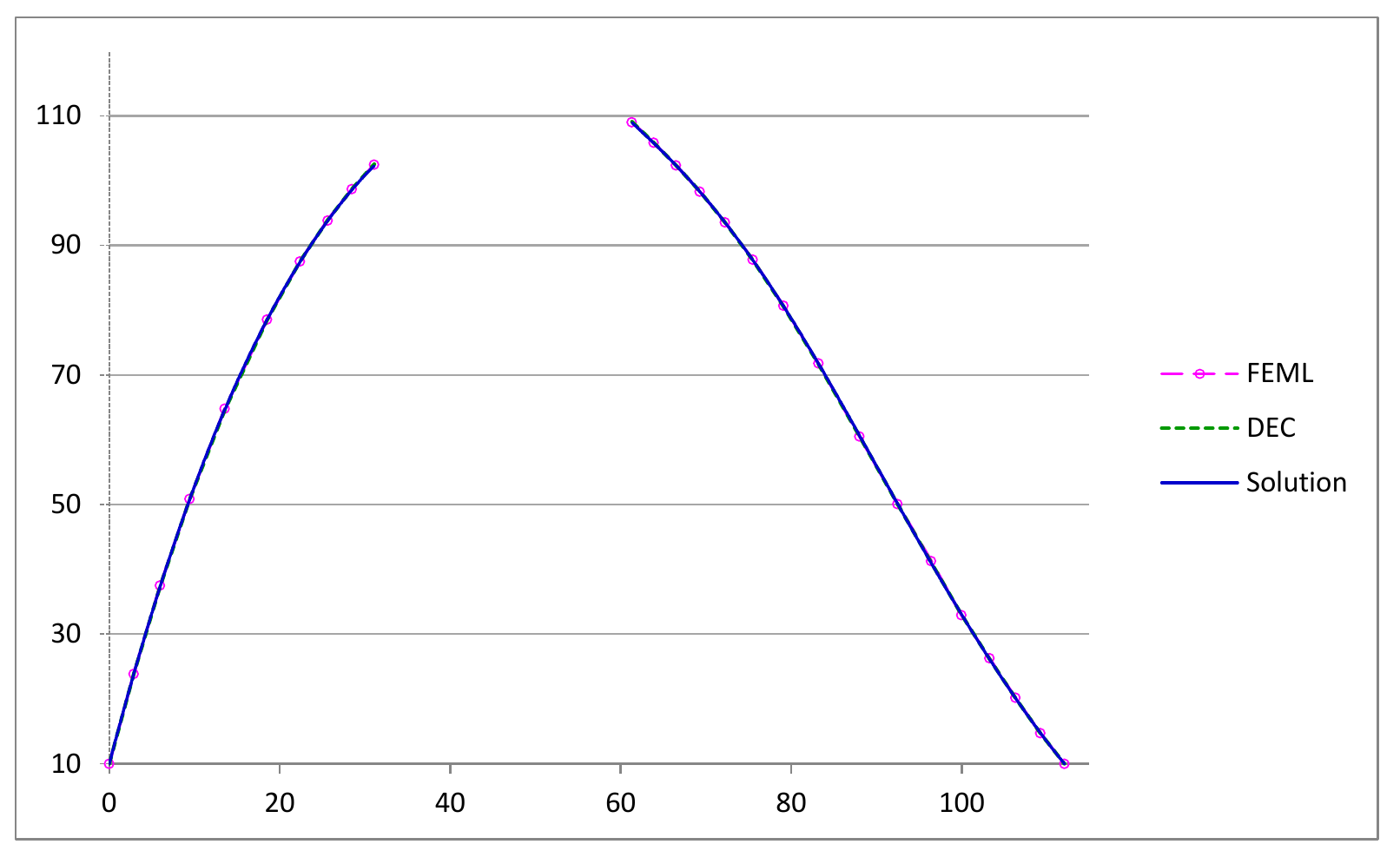}%
    \hspace{.1in}
    \includegraphics[width=0.4\textwidth]{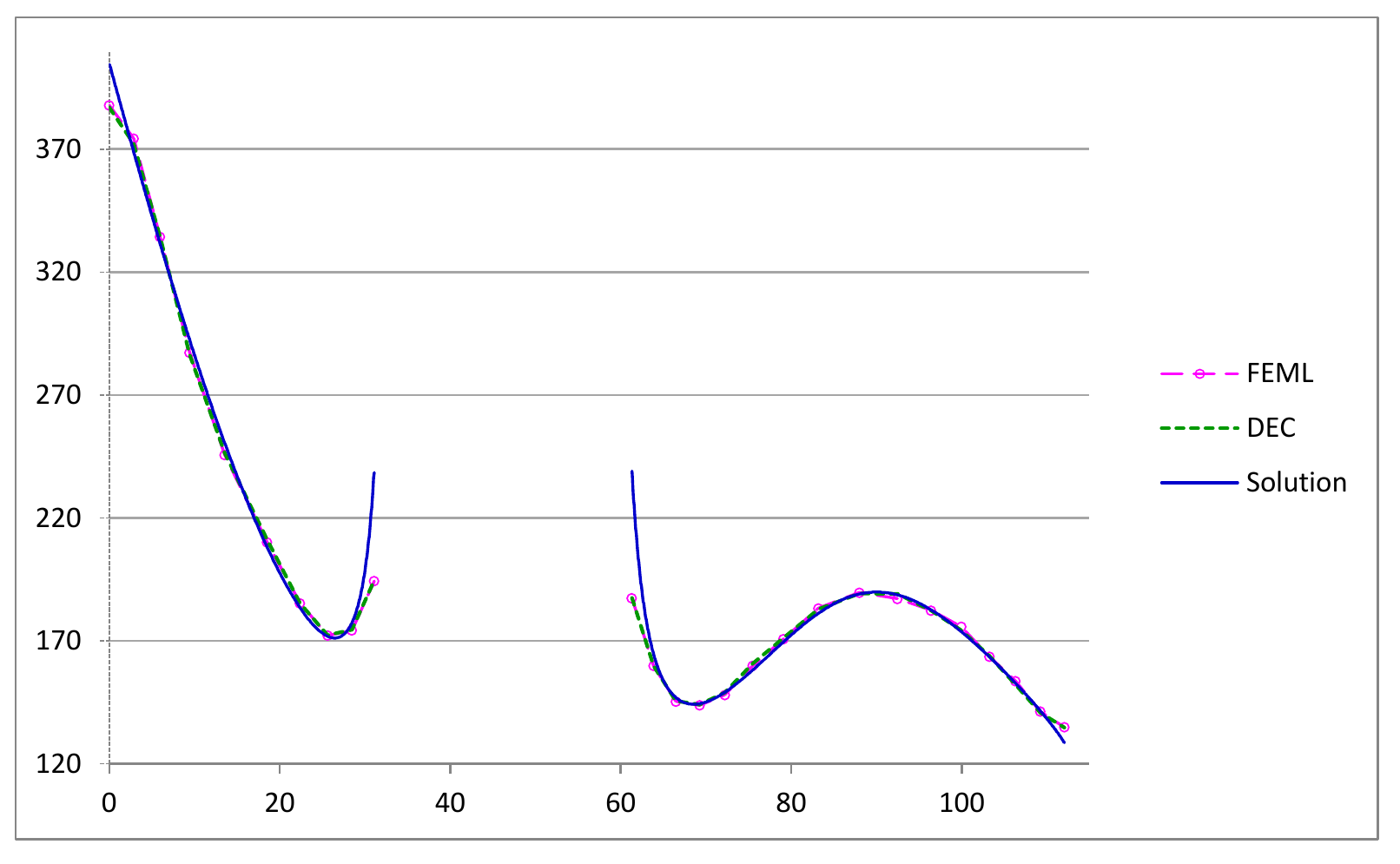}%
    \captionof{figure}{Temperature and Flux magnitude graphs along a horizontal line crossing the elliptical boundary for the first two meshes.}\label{fig: third example diametral graphs}
    \end{center}
As can be seen from Table \ref{table: third example numerical results} and Figure \ref{fig: third example diametral graphs}, 
the performance of DEC is very similar to that of FEML in this example. 
As expected, when the mesh is refined, the two methods converge to the same values.

\section{Conclusions}\label{sec: conclusions}

DEC is a relatively recent discretization scheme for PDE which takes into account the geometric and analytic 
features of the operators and the domains involved. 
The main contributions of this paper are the following:
\begin{enumerate}
 \item We have presented 2D DEC in a simplified manner, 
avoiding references to the theory of differential forms and motivating geometrically the new operators.
 \item We have carried out a numerical comparison between DEC and FEML by solving the 2D Poisson equation
 on two cirved domains. The numerical experiments show the solutions obtained with DEC on coarse meshes are as good or better as those of FEML.
 On the other hand, the experiements also show numerical convergence.
\item The computational cost of DEC is similar to that of FEML.
\end{enumerate}

\bigskip

{\em Acknowledgements}. The second named author was partially supported by a grant of CONACYT, and
would like to thank the International Centre for Numerical Methods in Engineering (CIMNE) and the University of Swansea
for their hospitality. We gratefully acknowledge the support of NVIDIA Corporation with the donation of the Titan X Pascal GPU used for this research.

\bigskip

{\small
\renewcommand{\baselinestretch}{0.5}
 }

\end{document}